\documentclass[12pt,twoside]{amsart}
\usepackage{amssymb,amsmath,amsfonts,amsthm}
\usepackage{mathrsfs}
\usepackage[X2,T1]{fontenc}
\usepackage[applemac]{inputenc}
\usepackage{cite}
\usepackage{latexsym}
\usepackage{verbatim}
\usepackage{soul}
\usepackage{ulem}
\usepackage[dvipsnames]{xcolor}
\def\Im{\mathop{\rm Im}\nolimits}
\def\Re{\mathop{\rm Re}\nolimits}

\def\R{\mathbb R}

\def\N{\mathbb N}

\def\ds{\displaystyle}

\newcommand\dslash{d\llap {\raisebox{.9ex}{$\scriptstyle-\!$}}}
\usepackage{colortbl}

\newcommand{\beqsn}{\arraycolsep1.5pt\begin{eqnarray*}}
\newcommand{\eeqsn}{\end{eqnarray*}\arraycolsep5pt}
\newcommand{\beqs}{\arraycolsep1.5pt\begin{eqnarray}}
\newcommand{\eeqs}{\end{eqnarray}\arraycolsep5pt}

\newtheorem{theorem}{Theorem}[section]
\newtheorem{lemma}[theorem]{Lemma}

\newtheorem{proposition}[theorem]{Proposition}
\newtheorem{definition}[theorem]{Definition}

\newtheorem{remark}[theorem]{Remark}

\catcode`\@=11

\renewcommand{\section}%
   {\setcounter{equation}{0}\@startsection {section}{1}{\z@}{-3.5ex plus -1ex
  minus -.2ex}{2.3ex plus .2ex}{\Large\bf}}

\title[ ]{Degenerate 3-evolution equations in Gevrey classes}

\author[A. Arias Junior]{Alexandre Arias Junior}
\address{Department of Computing and Mathematics \\ Universidade de S\~ao Paulo\\
	Ribeir\~ao Preto\\
	Brazil}
\email{alexandre.ariasjunior@usp.br}

\author[A. Ascanelli]{Alessia Ascanelli}
\address{Dipartimento di Matematica ed Informatica\\Universit\`a di Ferrara\\
Via Machiavelli 30\\
44121 Ferrara\\
Italy}
\email{alessia.ascanelli@unife.it}

\begin{document}


\begin{abstract}
	We consider the Cauchy problem for third-order evolution differential operators with variable coefficients, depending on time $t\in [0,T]$ and space $x\in\R$, where the leading coefficient 
$a_3(t)$ vanishes at $t = 0$ with finite order. We establish sufficient conditions on the behavior of the lower order coefficients $a_j(t,x)$ $j=1,2$ as $t \to 0^{+}$ and $|x| \to \infty$ that ensure well-posedness in $L^2(\R)$, $H^{\infty}(\R)$ and Gevrey-type spaces.
 \end{abstract}

\maketitle

\section{Introduction and state of the art}
In this paper we deal with the Cauchy problem 
\beqs
\label{genCP}
\begin{cases}
	P(t,x,D_t,D_x)u(t,x)=f(t,x) & (t,x)\in[0,T]\times\R\cr
	u(0,x)=g(x) & x\in\R
\end{cases}
\eeqs
for the operator
\beqs\label{P}
P = D_t + a_3(t)D^{3}_{x} + a_2(t,x) D^{2}_{x} + a_1(t,x)D_x + a_0(t,x),	
\eeqs
with $a_3\in C([0,T],\R)$ $a_{j} \in C([0,T];\mathcal{B}^{\infty}(\R)), j=0,1,2$, and $\mathcal{B}^{\infty}(\R)$ is the space of complex-valued functions bounded on $\R$ together with all their derivatives. The operator $P$ is known in the literature as a {\it $3-$evolution operator with real characteristics}, cf. \cite{mizohata2014}; indeed, the principal symbol (in the Petrowski sense) $p_3(t,\tau,\xi)=\tau+a_3(t)\xi^3$ admits the real characteristic root $\tau=-a_3(t)\xi^3$. The assumption $a_3(t) \in \R$ agrees with the necessary condition $\Im a_3(t)\leq 0\ \forall t\in [0,T]$ for well-posedness in $H^\infty(\R)= \bigcap_{m \in \R} H^m(\R)$ given by \cite[Theorem 3, page 31]{mizohata2014}. The fact that the lower order terms $a_j$ are complex-valued makes the problem challenging.

\smallskip

In this paper we are mainly interested in the case when the leading term $a_3(t)$ of the equation may vanish at some $t\in[0,T]$, that is, the equation degenerates and the behavior of the lower order terms $a_j(t,x),\ j=1,2,$ determines the space where the Cauchy problem \eqref{genCP} is well-posed: as we will see, we can have no $H^\infty$ well-posedness because of a degeneration with respect to $t$, or because of the behaviour at $|x| \to \infty$ of the coefficients, or both. In such cases, we shall seek well-posedness results in the Gevrey-type space
$$
{H}^\infty_{\theta} (\R): = \bigcup_{\rho >0}H^m_{\rho; \theta}(\R),
$$
where, given $\theta \geq 1, m, \rho \in \R,$ $H^{m}_{\rho; \theta} (\R)$ denotes the Gevrey-Sobolev space
$$
H^{m}_{\rho; \theta} (\R) = \{ u \in \mathscr{S}'(\R) :  \langle D \rangle^{m} 
e^{\rho \langle D \rangle^{\frac{1}{\theta}}} u \in L^{2}(\R) \}
$$
and $\langle D \rangle^m$, $e^{\rho \langle D \rangle^{\frac{1}{\theta}}}$ denote the Fourier multipliers with symbols $\langle \xi \rangle^m$ and $e^{\rho \langle \xi \rangle^{\frac{1}{\theta}}}$ respectively. These spaces are Hilbert spaces with inner product
$$
\langle u, v \rangle_{H^{m}_{\rho;\theta}} = \, \langle\langle D \rangle^{m} e^{\rho\langle D \rangle^{\frac{1}{\theta}}}u, \langle D \rangle^{m} e^{\rho\langle D \rangle^{\frac{1}{\theta}}}v \rangle_{L^{2}}, \quad u, v \in H^{m}_{\rho;\theta}(\R).
$$
The relation between $H^{\infty}_{\theta}(\R)$ and the Gevrey classes is given by $ G_0^\theta(\R) \subset  {H}^\infty_{\theta} (\R) \subset G^\theta(\R),$ where $G^\theta(\R)$ denotes the space of all smooth functions $f$ on $\R$ such that 
\begin{equation}
\label{gevestimate}
\sup_{\alpha \in \N^n} \sup_{x \in \R} h^{-|\alpha|} \alpha!^{-\theta} |\partial^\alpha f(x)| < +\infty
\end{equation} 
for some $h >0$, and $G_0^\theta(\R)$ is the space of all compactly supported functions contained in $G^\theta(\R)$. Notice that $H^{m}_{\rho;\theta}(\R)$ recaptures the the Sobolev space $H^m(\R)$ when $\rho=0$.

\medskip

If $a_3(t)$ never vanishes, i.e., under the assumption $a_3(t)\geq C>0$  (or $a_3(t)\leq - C<0$) $\forall t\in [0,T]$, the $H^\infty$ theory for the class of operators \eqref{P} is well-estabilished: necessary and sufficient conditions for $H^\infty$ well-posedness have been given, respectively, in  \cite{ascanelli_chiara_zanghirati_necessary_condition_for_H_infty_well_posedness_of_p_evo_equations} (where a more general class of $p$-evolution equations is considered),\cite{CC} and \cite{ascanelli_chiara_zhanghirati_well_posedness_of_cauchy_problem_for_p_evo_equations}. 
Roughly speaking, these conditions state that if the coefficient $\Im a_2$ decays as $x\to\infty$ as $\langle x\rangle^{-\sigma_2}$ with $\sigma_2\geq 1$ and $\Im a_1$ decays like $\langle x\rangle^{-\sigma_1}$ with $\sigma_1\geq 1/2$, then under some technical assumptions the Cauchy problem is well-posed in $H^\infty$; in the case of slowly decaying coefficients, i.e. $\sigma_2\in (0,1)$ or $\sigma_1\in (0,1/2)$, there's no hope for $H^\infty$ well-posedness.

In this second situation, necessary and sufficient decay conditions for well-posedness in $H^\infty_\theta(\R)$ with $\theta > 1$ of \eqref{genCP} have been recently stated, again under the assumption $a(t)_3\geq C>0$ $\forall t\in [0,T]$. On the one hand, from \cite{AACpevolGevreynec} (where the case of a more general class of $p-$evolution equations is considered), we know that if 
$$
		\Im \, a_{j}(t,x) \geq A \langle x \rangle^{-\sigma_{j}}, \quad x > R \, (\textrm{or} \, \, \, x <-R), \, t \in [0,T], \, j=1,2, 
$$
for some $R, A>0$ and $\sigma_{j} \in [0,1]$, 
$$
		|\partial^{\beta}_{x}a_{j}(t,x)| \leq C^{\beta + 1} \beta! \langle x \rangle^{-\beta}, \quad x \in \R, \, t \in [0,T], \,\, \beta \in \N_0, \quad j = 1,2,
$$
and the Cauchy problem \eqref{genCP} is well-posed in $
{H}^{\infty}_{\theta}(\R) $ for some $\theta > 1$, then
	\begin{equation}\nonumber
		\Xi := \max \{2(1-\sigma_{2}), 2(1-\sigma_{1})-1\} \leq \frac{1}{\theta}.
	\end{equation}
Testing this condition on the model operator $D_t + D^3_x + i\langle x \rangle^{-\sigma_2}D^{2}_{x}$, $\sigma_2 \in (0,1)$, we immediately realize that there is no well-posedness in ${H}^\infty_\theta(\R)$ of the Cauchy problem when $\sigma_2 \leq \frac{1}{2}$. On the other hand, see \cite{AAC3evolGevrey} (and \cite{eliakim} for a generalization), if there exist $\theta_0>1$, $\sigma \in \left(\frac{1}{2},1\right)$, $C_2,C_1>0$, such that for every $ \beta \in \N_0,$ $(t,x) \in [0,T] \times \R$
$$
|\partial_x^\beta a_{2}(t,x)| \leq C_{{2}}^{\beta+1} \beta!^{\theta_0} \langle x\rangle^{-\sigma-\beta},\quad  |\partial_x^\beta a_{1}(t,x)| \leq C_{{1}}^{\beta+1} \beta!^{\theta_0} \langle x \rangle^{-\frac{\sigma}{2}-\beta}
$$
then the Cauchy problem  \eqref{genCP} is well-posed in ${H}_\theta^\infty(\R)$ for $\theta\in\left[\theta_0,\frac{1}{2(1-\sigma)}\right)$. 

\bigskip

Coming now to the case where $a_3(t)=0$ for some $t\in[0,T]$, what happens? As far as the authors know, there is no answer in literature to this question.
The problem of a vanishing leading term in $p$-evolution equations has been considered at the moment only in the case of Schr\"odinger type equations, that is,
$$
P=D_t+a_2(t) D_x^2+ia_1(t,x)D_x+ a_0(t,x),
$$
where also the $n$-dimensional case was treated. In \cite{CRJEECT} the authors assume that there exist positive constants $\ell>0$ and $0< k\leq \ell,$ $\sigma, c,C>0$ such that
\beqsn
ct^\ell\leq a_2(t)&\leq& Ct^\ell,\quad \forall t\in [0,T],
\\
|\Im a_1(t,x)|&\leq& Ct^k\langle x\rangle^{-\sigma},\quad  \forall (t,x)\in [0,T]\times\R,
\eeqsn
 and they show that the corresponding Cauchy problem is well-posed in $H^\infty$ if $\ell=k$ and $\sigma\geq 1,$ while if $k<\ell$ or/and $\sigma<1$ (and $a_1$ is Gevrey regular of index $\theta_0$) it is well-posed in $H^\infty_\theta$ for $\theta\in[\theta_0,\theta_{\ell,k,\sigma})$ with upper bound $\theta_{\ell,k,\sigma}$ depending on the interplay between the degeneracy in time ($\ell,k$) and the decay of $Im\,a_1$ in space ($\sigma$). We also remark that in \cite{FS} the atuhors studied smoothing estimates for degenerate in time Schr\"odinger operators (there, no well-posedness results in Gevrey-type spaces were achieved). 

\section{Main result}

Let us consider the Cauchy problem \eqref{genCP}, \eqref{P} where we assume that all the coefficients are continuous with respect to time and satisfy, for some positive constants $\ell,k,k',>0$, $\sigma_1,\sigma_2>0$, and for every $\beta\in\N$, $(t,x)\in[0,T] \times\R$,
\beqs\label{a3}
c_{a_3}t^{\ell}\leq |a_3(t)|&\leq &C_{a_3}t^{\ell}, \quad 0<c_{a_3}\leq C_{a_3},
\eeqs
\beqs
\label{a2}
|\partial^{\beta}_{x} a_2(t,x)| &\leq& t^{k} C^{\beta+1}_{a_2}\beta! \langle x \rangle^{-\sigma_2}, 
\\\label{a1}
|\partial^{\beta}_{x} a_1(t,x)| &\leq& t^{k'} C^{\beta+1}_{a_1} \beta! \langle x \rangle^{-\sigma_1}, 
\\\label{a0}
|\partial^{\beta}_{x} a_0(t,x)| &\leq& C^{\beta+1}_{a_0}\beta!,
\eeqs
that is, the coefficient $a_3(t)$ vanishes of finite order $\ell$ at $t = 0$ and \eqref{a2}, \eqref{a1} are the proposed Levi-type conditions to control the degeneracy of $a_3(t)$. 

We are mainly interested in the "effective" degeneracy in time, that is
$$
\ell > \min\{k, k'\},
$$
which means that as $t\to 0^+$ we have $t^\ell< \min\{t^k, t^{k'}\}$, i.e. the degenerate behavior in $t$ of the lower order terms is worse with respect to the degenerate behavior in $t$ of the leading term. However, we are going to state general results including also the complementary case which is more simple to be studied, as we will see.

Before stating the main result of this paper, we define the two following numbers depending on the vanishing order $\ell$, the Levi-type conditions $k, k'$ and the decay rates $\sigma_1, \sigma_2$:
\beqs\label{qq2}
q_2 = q_2(\ell, k, \sigma_2) =
\begin{cases}
 2\frac{(\ell-k)\sigma_2 + (k+1)(1-\sigma_2)}{\sigma_2(\ell-k)+(k+1)}
= 2 \left( 1 - \frac{\sigma_2(k+1)}{\sigma_2(\ell-k)+(k+1)}\right), \quad \sigma_2 \in (\frac12,1), \\
\frac{2(\ell-k)}{\ell+1}, \qquad\qquad\qquad\qquad\qquad\qquad\qquad\qquad\quad \sigma_2 \geq 1,
\end{cases}
\eeqs
\beqs\label{qq1}
q_1 = q_1(\ell, k', \sigma_1) = 
\begin{cases}
\frac{(k'+1)(1-2\sigma_1)+\sigma_1(\ell-k')}{\sigma_1(\ell-k')+(k'+1)} = 1 - \frac{2\sigma_1(k'+1)}{\sigma_1(\ell-k') + (k'+1)}, \quad \sigma_1 \in (0,1), \\
\frac{\ell-2k'-1}{\ell+1}, \qquad\qquad\qquad\qquad\qquad\qquad\qquad\quad\quad\, \sigma_1 \geq 1.
\end{cases}
\eeqs
We then set $q = \max\{q_2,q_1\}$. 

\begin{remark}\label{rem1}
	Consider $\ell > k$. Then $q_2 > 0$ for all $\sigma_2 > 0$. Moreover, if $\ell < 2k+1$ then 
	$$
	q_2 < 1 \iff \sigma_2>\frac{k+1}{2+3k-\ell}. 
	$$
	We also point out that in the case $\ell \geq 2k+1$ we have $q_2 \geq 1$ for all $\sigma_2 > 0$. Regarding the index $q_1$ we have $q_1 < 1$ for all $\sigma_1 > 0$. If $q_1 \leq 0$ we shall consider $q_1 = 0$ and we will understand $\frac{1}{q_1}$ as $\infty$.
\end{remark}

The main result of this paper, that is the upcoming Theorem \ref{main}, states that, if $a_2$ presents a degeneracy in time with degeneracy gap $\ell-k>0$, the Cauchy probem is well-posed in suitable $H^{\infty}_{\theta}(\R)$ spaces, and the time and space behaviors of the coefficients play together to determine the Gevrey indices $\theta$ where well-posedness holds.

\begin{theorem}\label{main}
Consider the Cauchy problem \eqref{genCP} for the operator \eqref{P} under the assumptions \eqref{a3}, \eqref{a2}, \eqref{a1}, \eqref{a0}. Assume that 
\begin{equation}\label{condition.ell.k}
\ell > k, \quad \ell < 2k+1,
\end{equation}
and 
\begin{equation}\label{lowbound}
	\sigma_2 > \frac{k+1}{2+3k-\ell}.
\end{equation}
Then the Cauchy problem \eqref{genCP} is well-posed in $H^{\infty}_{\theta}(\R)$ for all $\theta$ satisfying
$$
1 < \theta < \frac{1}{q} = \min \left\{ \frac{1}{q_2}, \frac{1}{q_1}\right\}.
$$
\end{theorem}

Now we consider the case $k \geq \ell > k'$, that is the case when only $a_1$ presents a degeneration in time with degeneracy gap $\ell-k'>0$. In this situation, we can prove well-posedness in Sobolev spaces if the coefficient $a_1$ presents a not-too-slow decay as $x\to \infty$ and the degeneracy gap $\ell-k'$ is small enough (see Remark \ref{remark_q_1_greater_zero}); otherwise, we still find $H^\infty_\theta$ well-posedness with upper bound for $\theta$ depending on $\sigma_2,\sigma_1,\ell,k'$.  More precisely:

\begin{theorem}\label{main2}
Consider the Cauchy problem \eqref{genCP} for the operator \eqref{P} under the assumptions \eqref{a3}, \eqref{a2}, \eqref{a1}, \eqref{a0} and suppose that only $a_1$ presents a degeneracy gap in time, that is
\beqs\label{29}
k \geq \ell > k'.
\eeqs
If $\sigma_2 \in (\frac{1}{2}, 1)$, then the Cauchy problem \eqref{genCP} is well-posed in $H^{\infty}_{\theta}(\R)$ for all $\theta$ satisfying 
$$
1 < \theta < \frac1q=\min \left\{ \frac{1}{2(1-\sigma_2)}, \frac{1}{q_1} \right\}.
$$
If $\sigma_2 \geq 1$ and $q_1 > 0$, then the Cauchy problem \eqref{genCP} is well-posed in $H^{\infty}_{\theta}(\R)$ for all $\theta$ satisfying 
$$
1 < \theta < \frac{1}{q_1}.
$$
Finally, if $\sigma_2 \geq 1$ and $q_1 \leq 0$ we have well-posedness in 
$$
\begin{cases}
	H^{\infty}(\R), \quad \sigma_2 = 1, \\
	L^2(\R), \quad \sigma_2 > 1.
\end{cases}
$$
\end{theorem}

\begin{remark}
Notice that under the assumptions of Theorem \ref{main2}, that is when no degeneracy gap at level 2 appears ($\ell\leq k$), condition \eqref{lowbound} is replaced by $\sigma_2 > \frac{1}{2}$, which has been proved to be a necessary condition for $H^\infty_\theta$ well-posedness, see \cite{AACpevolGevreynec}.
\end{remark}

\begin{remark}\label{remark_q_1_greater_zero}
	Under the assumption \eqref{29} of Theorem \ref{main2}
	$$
	q_1 > 0 \iff \sigma_1(\ell-k') + (1-2\sigma_1)(k'+1) = \sigma_1(\ell-3k'-2)  + (k'+ 1) > 0.
	$$
	Therefore, $q_1 \leq 0$ is only possible when $\ell < 3k' + 2$ and $\sigma_1 > \frac{1}{2}$. Hence if $\ell < 3k' + 2$ (the degeneracy gap $\ell - k' > 0$ is not too large), $\sigma_1 > \frac{1}{2}$ ($a_1$ has "not too slow" decay) and the number $q_1 \leq 0$, then we can get $H^{\infty}(\R)$ well-posedness for \eqref{genCP} (provided that the decay rate of the coefficient $a_2$ is $\sigma_2 \geq 1$) even with the degeneracy gap $\ell - k' > 0$. Observe that the when $\ell > k$ we cannot obtain well-posedness in $H^{\infty}(\R)$, for all $\sigma_2 > 0$. In this sense, the degeneracy gap $\ell-k > 0$ is stronger than $\ell-k'>0$.
\end{remark}

Let us finally consider the case with no degeneracy gaps in time: $k, k' \geq \ell$. In this situation, only the decay rates at infinity of $a_2,a_1$ determine the space of well posedenss for the Cauchy problem \eqref{genCP}; if the decays are not-too-slow, then Sobolev well-posedness holds. More precisely:

\begin{theorem}\label{main3} 
	Consider the Cauchy problem \eqref{genCP} for the operator \eqref{P} under the assumptions \eqref{a3}, \eqref{a2}, \eqref{a1}, \eqref{a0} and suppose that 
	\beqs
	\ell\leq \min\{k,k'\},
	\eeqs
	that is there is no degeneracy gap in time. Then
	\begin{itemize}
		\item if $\sigma_2 \in (\frac12, 1), \sigma_1 \in (0, \frac{1}{2})$, then  the Cauchy problem \eqref{genCP} is well-posed in $H^{\infty}_{\theta}(\R)$ for every $1<\theta < s$ with 
		$$
		s:= \min\left\{\frac{1}{2(1-\sigma_2)}, \frac{1}{1-2\sigma_1} \right\};
		$$
		\item  if $\sigma_2 \in (\frac12, 1), \sigma_1 \geq \frac{1}{2}$, then the Cauchy problem \eqref{genCP} is well-posed in $H^{\infty}_{\theta}(\R)$ for every $1<\theta < s$ with 
		$$
		s:=\frac{1}{2(1-\sigma_2)};
		$$
		\item if $\sigma_2 \geq 1, \sigma_1 \in (0, \frac{1}{2})$, then the Cauchy problem \eqref{genCP} is well-posed in $H^{\infty}_{\theta}(\R)$ for every $1<\theta < s$ with 
		$$
		s:= \frac{1}{1-2\sigma_1};
		$$
		\item if $\sigma_2 = 1, \sigma_1 \geq \frac{1}{2}$, then the Cauchy problem \eqref{genCP} is well-posed in $H^{\infty}(\R);$
		\item if $\sigma_2 > 1, \sigma_1 \geq \frac{1}{2}$, then the Cauchy problem \eqref{genCP} is well-posed in $L^{2}(\R).$
	\end{itemize}
\end{theorem}

\begin{remark}
	Notice that in \eqref{a2}, \eqref{a1}, \eqref{a0} we have considered analytic coefficients with respect to space only for simplicity's sake. Indeed, we can assume Gevrey coefficients of index $\theta_0>1$, i.e. we can ask
	$$
	|\partial^{\beta}_{x} a_2(t,x)| \leq t^{k} C^{\beta+1}_{a_2}\beta!^{\theta_0} \langle x \rangle^{-\sigma_2}, \quad
	|\partial^{\beta}_{x} a_2(t,x)| \leq t^{k'} C^{\beta+1}_{a_1}\beta!^{\theta_0} \langle x \rangle^{-\sigma_1}, \quad
	|\partial^{\beta}_{x} a_0(t,x)| \leq C^{\beta+1}_{a_0}\beta!^{\theta_0},
	$$
	and obtain the same results of Theorems \ref{main}, \ref{main2} and  \ref{main3} with lower bound of the Gevrey index given by $\theta > \theta_0$ instead of $\theta > 1$.
\end{remark}

\begin{remark}
The first thesis of Theorem \ref{main3} provides  a generalization of \cite{AAC3evolGevrey} (see also Remark 6.1 in \cite{eliakim}). There, the case $\sigma_1=\sigma_2/2$ has been considered and the authors obtain for $\theta$ the upper bound $1/(2(1-\sigma_2))$ which is clearly consistent with Theorem \ref{main3}. In the present paper we have so removed the link between $\sigma_2$ and $\sigma_1$ that are here completely independent. The last two thesis of Theorem \ref{main3} are already known, see \cite{CC, ascanelli_chiara_zhanghirati_well_posedness_of_cauchy_problem_for_p_evo_equations}.
\end{remark}

\begin{remark}
	It is interesting to test wether the constraints on the parameters $\ell, k, k', \sigma_2, \sigma_1$ in theorems \ref{main} and \ref{main2} are sharp or not. For instance, when $\ell \geq 2k+1$ (the degeneracy gap $\ell - k >0$ is large) we cannot conclude well-posedness in $L^2$, $H^{\infty}$ and $H^{\infty}_{\theta}$ for all $\theta > 1$, even assuming a very fast decay of $a_2$ and $a_1$ as $x\to \infty$. We observe that if $k \geq \ell$ (no degeneracy gap between the leading term and the part of order $2$) $H^\infty_\theta$ well-posendess is always possibile for some $\theta$, also in the case of a very large gap $\ell - k' > 0$ at level 1. We may observe an analogous situation for the Schr\"odinger case, see \cite{CRJEECT}. Therefore it is interesting to test wether $\ell < 2k+1$ is a technical condition or not in Theorem \ref{main}. Besides, it is also interesting to test if we have sharp spaces of well-posedness in theorems \ref{main} and \ref{main2}. We intend to study this topics in a future paper.
\end{remark}

	In the following sections we are going to prove first Theorem \ref{main}, which is the most challenging one. Indeed, in Section 4 we will start with the proof of Theorem \ref{main} in the most difficult case, that is when $\sigma_2,\sigma_1<1$, namely, when both degeneracy in time and slow decay in space of the coefficients occur. These two destroy well-posedness in Soboev spaces and both determine the class of Gevrey-type spaces in which well-posedness can be found. After that, still in Section 4, we consider the case when at least one between $\sigma_2,\sigma_1$ is greater than or equal to $1.$ Finally, in Section 5, we shall discuss the needed changes in the proof of Theorem \ref{main} to get Theorems \ref{main2} and \ref{main3}. The upcoming Section 3 contains some preliminary results concerning pseudodifferential calculus for operators of finite and infinite order that will be used throughout the paper.

\section{Preliminaries}\label{preliminaries}

This section is devoted to give some definitions and results that we are going to use throughout the text. We remark that all the results stated in this section are valid in $\R^n$, however we state them with $n=1$ for simplicity sake.

\subsection{Pseudodifferential operators of finite order}

As usual, for any $m \in \R$, $S^{m}(\R^2)$ stands for the space of all functions $p \in C^{\infty}(\R^2)$ such that for any $\alpha, \beta \in \N_0$ there exists $C_{\alpha,\beta} > 0$ such that 
\begin{equation}\label{eq_def_symbol}
|\partial^{\alpha}_{\xi} \partial^{\beta}_{x} p(x,\xi)| \leq C_{\alpha,\beta} \langle \xi \rangle^{m-\alpha}.
\end{equation}
The topology in $S^{m}(\R^2)$ is induced by the following family of seminorms
$$
|p|^{(m)}_\ell := \max_{\alpha \leq \ell, \beta \leq \ell} \sup_{x, \xi \in \R} |\partial^{\alpha}_{\xi} \partial^{\beta}_{x} p(x,\xi)| \langle \xi \rangle^{-m+|\alpha|}, \quad p \in S^{m}(\R), \, \ell \in \N_0.
$$
We associate to every symbol $p \in S^{m}(\R^2)$ the continuous operator $p(x,D): \mathscr{S}(\R) \to \mathscr{S}(\R)$ (Schwartz space of rapidly decreasing functions), known as pseudodifferential operator, defined by 
$$
p(x,D) u(x) = \int e^{i\xi x} p(x,\xi) \widehat{u}(\xi) \dslash\xi, \quad u \in \mathscr{S}(\R),
$$
where $\hat u$ denotes the Fourier transform of $u$ and $\dslash\xi=(2\pi)^{-1}d\xi$. We will often denote the operator $p(x,D)$ by $op(p(x,\xi))$ referring to its symbol. The next result gives the action of operators coming from symbols $S^{m}(\R^2)$ in the standard sobolev spaces $H^{s}(\R)$, $s \in \R$. For a proof we address the reader to Theorem $1.6$ on page $224$ of \cite{KG}.

\begin{theorem}\label{theorem_Calderon_Vaillancourt}[Calder\'on-Vaillancourt]
	Let $p \in S^{m}(\R^2)$. Then for any real number $s \in \R$ there exist $\ell := \ell(s,m) \in \N_0$ and $C:= C_{s,m} > 0$ such that 
	\begin{equation*}
		\| p(x,D)u \|_{H^{s}(\R)} \leq C |p|^{(m)}_{\ell} \| u \|_{H^{s+m}(\R)}, \quad \forall \, u \in H^{s+m}(\R).
	\end{equation*}
	Besides, when $m = s = 0$ we can replace $|p|^{(m)}_{\ell}$ by
	\begin{equation*}
		\max_{\alpha, \beta \leq 2} \sup_{x, \xi \in \R} |\partial^{\alpha}_{\xi} \partial^{\beta}_{x} p(x,\xi)|.
	\end{equation*}
\end{theorem}

Now we consider the algebra properties of $S^{m}(\R^2)$ with respect the composition of operators. Let $p_j \in S^{m_j}(\R^2)$, $j = 1, 2$, and define 
\begin{align}\label{eq_symbol_of_composition}
	q(x,\xi) &= Os- \iint e^{-iy \eta} p_1(x,\xi+\eta)p_2(x+y,\xi) dy \dslash\eta \\
	&= \lim_{\varepsilon \to 0} \iint e^{-iy\eta} p_1(x,\xi+\eta)p_2(x+y,\xi) e^{-\varepsilon^2y^2} e^{-\varepsilon^2\eta^2} dy \dslash\eta. \nonumber
\end{align} 

Then we have the following theorem (for a proof see Lemma $2.4$ on page $69$ and Theorem $1.4$ on page $223$ of \cite{KG}).

\begin{theorem}
	Let $p_j \in S^{m_j}(\R)$, $j = 1, 2$, and consider $q$ defined by \eqref{eq_symbol_of_composition}. Then $q \in S^{m_1+m_2}(\R^2)$ and $q(x,D) = p_1(x,D)\circ p_2(x,D)$. Moreover, the symbol $q$ has the following asymptotic expansion
	\begin{align}\label{eq_asymptotic_expansion_formula}
		q(x,\xi) = \sum_{\alpha < N} \frac{1}{\alpha!} \partial^{\alpha}_{\xi}p_{1}(x,\xi)D^{\alpha}_{x}p_2(x,\xi) + r_N(x,\xi), 
	\end{align} 
	where 
	$$
	r_N(x,\xi) = N\int_{0}^{1} \frac{(1-\theta)^{N-1}}{N!} \, Os - \iint e^{-iy\eta} \partial^{N}_{\xi} p_1(x,\xi+\theta\eta) D^{N}_{x} p_2(x+y,\xi)  dy\dslash\eta \, d\theta,
	$$
	and the seminorms of $r_N$ may be estimated in the following way: for any $\ell_{0} \in \N_0$ there exists $\ell_{1} := \ell_1(\ell_{0}) \in \N_0$ such that 
	$$
	|r_N|^{(m_1+m_2)}_{\ell_0} \leq C_{\ell_{0}} |\partial^{N}_{\xi}p_1|^{(m_1)}_{\ell_{1}} |\partial^{N}_{x}p_2|^{(m_2)}_{\ell_{1}}.
	$$
\end{theorem} 

We also recall the so-called SG symbol classes. These classes take into account the algebraic growth with respect to the variable $x$. Being more precise, given $m_1, m_2 \in \R$ we say that $p \in SG^{m_1, m_2}(\R^2)$ when for all $\alpha, \beta \in \N_0$ there exists a constant $C_{\alpha, \beta} > 0$ such that 
$$
|\partial^{\alpha}_{\xi} \partial^{\beta}_{x} p(x,\xi)| \leq C_{\alpha,\beta} \langle \xi \rangle^{m_1-\alpha} \langle x \rangle^{m_2-\beta}.
$$

We end this subsection recalling the so-called sharp-G{\aa}rding and Fefferman-Phong inequalities. 

\begin{theorem}\label{theorem_SG_sharp_garding}[sharp-G{\aa}rding inequality]
	Let $p \in S^{m}(\R^{2})$. If $Re\, p(x,\xi) \geq 0$ for all $x,\xi \in \R$ then there exists a constant $C > 0$, depending on a finite number of the constants $C_{\alpha,\beta}$ in \eqref{eq_def_symbol}, such that 
	$$
	Re\, \langle p(x,D) u, u \rangle_{L^{2}} \geq -C \| u \|_{H^{\frac{m-1}{2}}}, \quad u \in \mathscr{S}(\R).
	$$
\end{theorem} 

\begin{theorem}\label{theorem_fefferman_phong}[Fefferman-Phong inequality]
	Let $p \in S^{m}(\R^{2})$ be a real-valued symbol. If $ p(x,\xi) \geq 0$ for all $x,\xi \in \R$ then there exists a constant $C > 0$, depending on a finite number of the constants $C_{\alpha,\beta}$ in \eqref{eq_def_symbol}, such that 
	$$
	Re\, \langle p(x,D) u, u \rangle_{L^{2}} \geq -C \| u \|_{H^{\frac{m-2}{2}}}, \quad u \in \mathscr{S}(\R).
	$$
\end{theorem}

\subsection{Pseudodifferential operators of infinite order }\label{section_pseudodifferential_operators_of_infinite_order}

In this subsection we will recall the main results concerning pseudodifferential operators of infinite order. First we define the Gelfand-Shilov spaces and then we will pass to the classes of symbols. For the proofs of the results here stated we refer to Chapter $2$ of \cite{arias_phd_thesis} and to Appendix A of \cite{ACJMPA}.

The Gelfand-Shilov spaces have been originally introduced in the book \cite{GS2}, making part of larger class called spaces of type $\mathcal{S}$. More to the point, given $\theta \geq 1$ and $A, B > 0$ we say that a smooth function $f$ belongs to $\mathcal{S}^{\theta, A}_{\theta, B} (\R)$ if there is a constant $C > 0$  such that 
$$
|x^{\beta} \partial^{\alpha}_{x}f(x)| \leq C A^{\alpha} B^{\beta}\alpha!^{\theta} \beta!^{\theta}, 
$$
for every $\alpha, \beta \in \N_{0}$ and $x \in \R$. The norm 
$$
\| f \|_{\theta, A, B} \,= \sup_{\overset{x \in \R}{\alpha, \beta \in \N^{n}_{0}}} |x^{\beta} \partial^{\alpha}_{x}f(x)|A^{-\alpha}B^{-\beta}\alpha!^{-\theta}\beta!^{-\theta} 
$$
turns $\mathcal{S}^{\theta, A}_{s, B}(\R)$ into a Banach space. We define 
$$
\mathcal{S}_{\theta}(\R) = \bigcup_{A,B >0} \mathcal{S}^{\theta, A}_{s, B} (\R)
$$ 
and we can equip it with the inductive limit topology of the Banach spaces $\mathcal{S}^{\theta, A}_{\theta, B}(\R)$. We also consider the projective version, that is 
$$
\Sigma_{\theta} (\R) = \bigcap_{A, B > 0} \mathcal{S}^{\theta, A}_{\theta, B} (\R)
$$
equipped with the projective limit topology. The following inclusions are continuous (for every $\varepsilon > 0$)
$$
\Sigma_{\theta} (\R) 
\subset \mathcal{S}_{\theta}(\R) \subset \Sigma_{\theta+\varepsilon} (\R).
$$
Concerning the action of the Fourier transform we have the following isomorphisms
$$
\mathcal{F}: \Sigma_{\theta} (\R) \to \Sigma_{\theta} (\R), \quad 
\mathcal{F}: \mathcal{S}_{\theta}(\R) \to \mathcal{S}_{\theta}(\R).
$$

Now we pass to the symbol classes of finite and infinite order with Gevrey regular symbols.

\begin{definition} Let $m, s \in \R$ and $\mu, \theta > 1$. Then
	\begin{itemize}
	\item[(i)] $ S^{m}_{\mu}(\R^{2})$ denotes the space of all functions $p \in C^{\infty}(\R^{2})$ for which there exist $C, A > 0$ such that 
	$$
	|\partial^{\alpha}_{\xi} \partial^{\beta}_{x} a(x,\xi)| \leq C A^{\alpha+\beta} \alpha!^{\mu} \beta!^{\mu} \langle \xi \rangle^{m-\alpha},
	$$
	for every $\alpha, \beta \in \N_{0}$ and $\xi, x \in \R$.
	
	\item[(ii)] $ SG^{m,s}_{\mu}(\R^{2})$ denotes the space of all functions $p \in C^{\infty}(\R^{2})$ for which there exist $C, A > 0$ such that 
	$$
	|\partial^{\alpha}_{\xi} \partial^{\beta}_{x} a(x,\xi)| \leq C A^{\alpha+\beta} \alpha!^{\mu} \beta!^{\mu} \langle \xi \rangle^{m-\alpha} \langle x \rangle^{s-\beta},
	$$
	for every $\alpha, \beta \in \N_{0}$ and $\xi, x \in \R$.
	
	\item [(iii)] $S^{\infty}_{\mu; \theta}(\R^2)$ denotes the space of all functions $p \in C^{\infty}(\R^{2})$ for which there exist $C, c, A > 0$ such that 
	$$
	|\partial^{\alpha}_{\xi} \partial^{\beta}_{x} a(x,\xi)| \leq C A^{\alpha+\beta} \alpha!^{\mu} \beta!^{\mu} \langle \xi \rangle^{-\alpha} e^{c|\xi|^{\frac{1}{\theta}}},
	$$
	for every $\alpha, \beta \in \N_{0}$ and $\xi, x \in \R$.
	
	\item [(iv)] $SG^{\infty,s}_{\mu; \theta}(\R^{2})$ denotes the space of all functions $p \in C^{\infty}(\R^{2})$ for which there exist $C, c, A > 0$ such that 
	$$
	|\partial^{\alpha}_{\xi} \partial^{\beta}_{x} a(x,\xi)| \leq C A^{\alpha+\beta} \alpha!^{\mu} \beta!^{\mu} \langle \xi \rangle^{-\alpha} \langle x \rangle^{s-\beta} e^{c|\xi|^{\frac{1}{\theta}}},
	$$
	for every $\alpha, \beta \in \N^{n}_{0}$ and $\xi, x \in \R$.
	\end{itemize}
\end{definition}

\begin{remark}
	We have the following inclusions $S^{m}_{\mu} \subset S^{\infty}_{\mu;\theta}$ and $SG^{m,s}_{\mu} \subset SG^{\infty,s}_{\mu;\theta}$. Moreover, if $s \leq 0$ then $SG^{m,s}_{\mu} \subset S^{m}_{\mu}$ and $SG^{\infty,s}_{\mu} \subset S^{\infty}_{\mu}$.
\end{remark}

We have the following continuity results.

\begin{proposition}
	Let $\theta > \mu \geq 1$ and $p$ either in $SG^{\infty, s}_{\mu; \theta}(\R^{2})$ or in $S^{\infty}_{\mu;\theta}(\R^n)$. Then the operator $p(x,D)$ is continuous on $\Sigma_{\theta}(\R)$ and it extends to a continuous map on $(\Sigma_{\theta})'(\R)$).
\end{proposition}

\begin{proposition}\label{theorem_continuity_gevrey_sobolev_finite_order}
	Let $p \in S^{m'}_{\mu}(\R^{2})$ for some $m' \in \R$. Then for every $m, \rho \in \R$ and $\theta$ such that $\theta > 2\mu - 1$ the operator $p(x, D)$ maps $H^{m}_{\rho; \theta}(\R)$ into $H^{m-m'}_{\rho; \theta}(\R)$ continuously.
\end{proposition}

In order to recall the pseudodifferential calculus for infinite order operators, we need the notion of asymptotic expansion. 

\begin{definition}\label{definition_classes_of_asymptotic_expansions} We say that: 
	\begin{itemize}
		\item [(i)] $\sum\limits_{j \geq 0} a_j \in FS^{\infty}_{\mu; \theta}$ if $a_j(x, \xi) \in C^{\infty}(\R^{2})$ and there are $C, c, B > 0$ satisfying
		$$
		|\partial^{\alpha}_{\xi}\partial^{\beta}_{x} a_j(x, \xi)| \leq C^{\alpha + \beta + 2j + 1} \alpha!^{\mu} \beta!^{\mu} j!^{2\mu -1} 
		\langle \xi \rangle^{- \alpha - j} e^{c|x|^{\frac{1}{\theta}}}, 
		$$
		for every $\alpha, \beta \in \N_{0}$, $j \geq 0$ and $\langle \xi \rangle \geq Bj^{2\mu - 1}$;
		
		\item [(ii)] $\sum\limits_{j \geq 0} a_j \in FSG^{\infty,s}_{\mu; \theta}$ if $a_j(x, \xi) \in C^{\infty}(\R^{2})$ and there are $C, c, B > 0$ satisfying
		$$
		|\partial^{\alpha}_{\xi}\partial^{\beta}_{x} a_j(x, \xi)| \leq C^{\alpha + \beta + 2j + 1} \alpha!^{\mu} \beta!^{\mu} j!^{2\mu -1} 
		\langle \xi \rangle^{- \alpha - j} \langle x \rangle^{s-\beta-j}e^{c|x|^{\frac{1}{\theta}}}, 
		$$
		for every $\alpha, \beta \in \N_{0}$, $j \geq 0$ and $\langle \xi \rangle \geq Bj^{2\mu - 1}$ or $\langle x \rangle \geq Bj^{2\mu - 1}$;
		
		\item [(iii)] $\sum\limits_{j \geq 0} a_j \in FS^{m}_{\mu}$ if $a_j(x, \xi) \in C^{\infty}(\R^{2})$ and there are $C, c, B > 0$ satisfying
		$$
		|\partial^{\alpha}_{\xi}\partial^{\beta}_{x} a_j(x, \xi)| \leq C^{\alpha + \beta + 2j + 1} \alpha!^{\mu} \beta!^{\mu} j!^{2\mu -1} 
		\langle \xi \rangle^{m - \alpha - j}, 
		$$
		for every $\alpha, \beta \in \N_{0}$, $j \geq 0$ and $\langle \xi \rangle \geq Bj^{2\mu - 1}$;
		
		\item [(iv)] $\sum\limits_{j \geq 0} a_j \in FSG^{m,s}_{\mu}$ if $a_j(x, \xi) \in C^{\infty}(\R^{2})$ and there are $C, c, B > 0$ satisfying
		$$
		|\partial^{\alpha}_{\xi}\partial^{\beta}_{x} a_j(x, \xi)| \leq C^{\alpha + \beta + 2j + 1} \alpha!^{\mu} \beta!^{\mu} j!^{2\mu -1} 
		\langle \xi \rangle^{m - \alpha - j} \langle x \rangle^{s-\beta-j}, 
		$$
		for every $\alpha, \beta \in \N_{0}$, $j \geq 0$ and $\langle \xi \rangle \geq Bj^{2\mu - 1}$  or $\langle x \rangle \geq Bj^{2\mu - 1}$.
	\end{itemize}	   
\end{definition}

\begin{definition}
	Let $\sum\limits_{j \geq 0} a_j$, $\sum\limits_{j \geq 0} b_j$ in $FS^{\infty}_{\mu; \theta}$. We say that $\sum\limits_{j \geq 0} a_j \sim \sum\limits_{j \geq 0} b_j$ in $FS^{\infty}_{\mu; \theta}$ if  there exist $C, c, B > 0$ satisfying 
	$$
	|\partial^{\alpha}_{\xi}\partial^{\beta}_{x} \sum_{j < N} (a_j - b_j) (x, \xi)| \leq C^{\alpha + \beta + 2N + 1} \alpha!^{\mu} \beta!^{\mu} N!^{2\mu - 1} 
	\langle \xi \rangle^{- |\alpha| - N} e^{c|\xi|^{\frac{1}{\theta}}}, 
	$$
	for every $\alpha, \beta \in \N^{n}_{0}$, $N \geq 1$ and $\langle \xi \rangle \geq B N^{2\mu-1}$. Analogous definition for the class $FS^{m}_{\mu, \nu}$.
\end{definition}

\begin{definition}
Let $\sum\limits_{j \geq 0} a_j$, $\sum\limits_{j \geq 0} b_j$ in $FSG^{\infty, s}_{\mu; \theta}$. We say that $\sum\limits_{j \geq 0} a_j \sim \sum\limits_{j \geq 0} b_j$ in $FSG^{\infty, s}_{\mu; \theta}$ if  there are $C, c, B > 0$ satisfying 
$$
|\partial^{\alpha}_{\xi}\partial^{\beta}_{x} \sum_{j < N} (a_j - b_j) (x, \xi)| \leq C^{\alpha + \beta + 2N + 1} \alpha!^{\mu} \beta!^{\mu} N!^{2\mu - 1} 
\langle x \rangle^{s - |\beta| - N} \langle \xi \rangle^{-|\alpha| - N} e^{c|\xi|^{\frac{1}{\theta}}}, 
$$
for every $\alpha, \beta \in \N_{0}$, $N \geq 1$ and $\langle \xi \rangle \geq B N^{2\mu-1}$ or $\langle x \rangle \geq B N^{2\mu-1}$. Analogous definition for the class $FSG^{m,s}_{\mu}$.
\end{definition}

\begin{remark}
	If $\sum\limits_{j \geq 0} a_j \in FS^{\infty}_{\mu; \theta}$, then $a_0 \in S^{\infty}_{\mu; \theta}$. Given $a \in S^{\infty}_{\mu; \theta}$ and setting $b_0 = a$, $b_j = 0$, $j \geq 1$, we have $a = \sum\limits_{j \geq 0}b_j$. Hence we can consider $S^{\infty}_{\mu; \theta}$ as a subset of $FS^{\infty}_{\mu; \theta}$. Analogous considerations for the other classes.
\end{remark} 

\begin{proposition}\label{Proposition_existence_of_a_stymbol_which_has_an_given_asymptotica_expansion}
	Given $\sum\limits_{j \geq 0} a_j \in FS^{\infty}_{\mu; \theta}$, there exists $a \in S^{\infty}_{\mu; \theta}$ such that $a \sim \sum\limits_{j \geq 0} a_j$ in $FS^{\infty}_{\mu; \theta}$. Analogous result for the other classes.
\end{proposition}
\begin{definition}
Let $r \geq 1$. Then $\mathcal{K}_{r}$ denotes the space of all symbols $p \in S^{\infty}_{\mu;\theta}$ such that there exist $C, c > 0$ such that
\begin{equation*}
	|\partial^{\alpha}_{\xi}\partial^{\beta}_{x}p(x,\xi)| \leq C^{\alpha+\beta+1}\alpha!^{r} \beta!^{r} e^{-c|\xi|^{\frac{1}{r}}}, \quad 
	x,\xi \in \R, \, \alpha,\beta \in \N_0.
\end{equation*}
\end{definition}
\begin{proposition}
	Let $a, b \in S^{\infty}_{\mu; \theta}$ and $\sum\limits_{j \geq 0} a_j \in FS^{\infty}_{\mu; \theta}$. If $a \sim \sum\limits_{j \geq 0} a_j$ in $FS^{\infty}_{\mu; \theta}$, $b \sim \sum\limits_{j \geq 0} a_j$ in $FS^{\infty}_{\mu; \theta}$ and $\theta > 2\mu-1$, then $a-b \in \mathcal{K}_{2\mu-1}$. 
	 Analogous result for the class $FS^{m}_{\mu}$.
\end{proposition}

\begin{proposition}
	Let $a, b \in SG^{\infty, s}_{\mu; \theta}$ and and $\sum\limits_{j \geq 0} a_j \in FSG^{\infty,s}_{\mu; \theta}$.
	If $a \sim \sum\limits_{j \geq 0} a_j$ in $FSG^{\infty,s}_{\mu; \theta}$, $b \sim \sum\limits_{j \geq 0} a_j$ in $FSG^{\infty,s}_{\mu; \theta}$ and $\theta > 2\mu-1$,
	then $a-b \in \mathcal{S}_{2\mu-1}$. Analogous result for the class $FSG^{m, s}_{\mu; \theta}$.
\end{proposition}

Concerning the symbolic calculus we have:

\begin{theorem}\label{theorem_symbolic_calculus_of_infinte_order}
	Let $p \in S^{\infty}_{\mu; \theta}(\R^{2})$, $q \in S^{\infty}_{\mu; \theta}(\R^{2})$ with $\theta > 2\mu - 1$. Then the $L^{2}$ adjoint $p^{*}$ and the composition $p\circ q$ have the following structure:
	\begin{itemize}	
		\item[-] $p^{*}(x,D) = a(x, D) + r(x,D)$ where $r\in \mathcal{K}_{2\mu-1}(\R^{2})$, $a \in S^{\infty}_{\mu; \theta}(\R^{2})$, and
		$$
		a(x,\xi) \sim \sum_{\alpha} \frac{1}{\alpha!} \overline{\partial^{\alpha}_{\xi}D^{\alpha}_xp(x,\xi)} \,\, \text{in} \,\,
		FS^{\infty}_{\mu; \theta}(\R^{2}).
		$$
		\item[-]$p(x,D)\circ q(x, D) = b(x,D) + s(x,D)$, where $s \in \mathcal{K}_{2\mu-1}(\R^{2})$, $b \in S^{\infty}_{\mu; \theta}(\R^{2})$ and
		$$
		b(x, \xi) \sim \sum_{\alpha} \frac{1}{\alpha!} \partial^{\alpha}_{\xi}p(x,\xi) D^{\alpha}_{x}q(x,\xi) \,\, \text{in} \,\,
		FS^{\infty}_{\mu; \theta}(\R^{2}).
		$$
	\end{itemize}
	Analogous result for the class $S^{m}_{\mu}(\R^{2})$.
\end{theorem}

\begin{remark}
	Let $\mu > 1$ and $\theta > 2\mu-1$, $r \in \mathcal{K}_{2\mu-1}$ and $p \in S^{\infty}_{\mu;\theta}(\R^2)$. Then $r^{*}(x,D)$, $p(x,D)\circ r(x,D)$ and $r(x,D) \circ p(x,D)$  are operators given by symbols in $\mathcal{K}_{2\mu-1}$. 
\end{remark}

\begin{theorem}\label{theorem_symbolic_SG_calculus_of_infinte_order}
	Let $p \in SG^{\infty, s}_{\mu; \theta}(\R^{2})$, $q \in SG^{\infty, s'}_{\mu;\theta}(\R^{2})$ with $\theta > 2\mu - 1$. Then the $L^{2}$ adjoint $p^{*}$ and the composition $p\circ q$ have the following structure:
	\begin{itemize}	
		\item[-] $p^{*}(x,D) = a(x, D) + r(x,D)$ where $r\in \mathcal{S}_{2\mu-1}(\R^{2})$, $a \in SG^{\infty, s}_{\mu; \theta}(\R^{2})$, and
		$$
		a(x,\xi) \sim \sum_{\alpha} \frac{1}{\alpha!} \partial^{\alpha}_{\xi}D^{\alpha}_x \overline{p(x,\xi)} \,\, \text{in} \,\,
		FSG^{\infty, s}_{\mu; \theta}(\R^{2});
		$$
		\item[-] $p(x,D)\circ q(x, D) = b(x,D) + w(x,D)$, where $w \in \mathcal{S}_{2\mu-1}(\R^{2})$, $b \in SG^{\infty,s+s'}_{\mu; \theta}(\R^{2})$ and
		$$
		b(x, \xi) \sim \sum_{\alpha} \frac{1}{\alpha!} \partial^{\alpha}_{\xi}p(x,\xi) D^{\alpha}_xq(x,\xi) \,\, \text{in} \,\,
		FSG^{\infty, s+s'}_{\mu; \theta}(\R^{2}).
		$$
	\end{itemize}
	Analogous result for the class $SG^{m,s}_{\mu}(\R^{2})$.
\end{theorem}

\begin{remark}
	If $\mu > 1$ and $\theta > 2\mu-1$. Let $r \in \mathcal{S}_{2\mu-1}$ and $p \in SG^{\infty,s}_{\mu;\theta}(\R^2)$. Then $r^{*}(x,D)$, $p(x,D)\circ r(x,D)$ and $r(x,D) \circ p(x,D)$  are operators given by symbols in $\mathcal{S}_{2\mu-1}$. 
\end{remark}

\subsection{Symbols of the form $e^{\lambda}$ and the reverse operator}

In the proof of Theorem \ref{main} we are going to consider operators of infinite order of the form $e^\lambda(x,D)$, where $\lambda$ is a Gevrey regular symbol. In this case, a simple application of the Fa\`a di Bruno's formula gives the following.

\begin{proposition}\label{proposition_exponential_of_symbols_of_finite_order}
	If $\lambda(x,\xi) \in S^{\frac{1}{\theta}}_{\mu}(\R^{2})$, then $e^{\pm\lambda(x,\xi)} \in S^{\infty}_{\mu; \theta}(\R^{2})$. If $\lambda(x,\xi) \in SG^{\frac{1}{\theta},0}_{\mu}(\R^{2})$, then $e^{\pm\lambda(x,\xi)} \in SG^{\infty,0}_{\mu; \theta}(\R^{2})$.
\end{proposition}

The so-called reverse operator $^{R}\{e^{\pm \lambda}\}(x,D)$ of $e^{\pm \lambda}(x,D)$ was introduced in \cite[Proposition 2.13]{KW} as the transposed of $e^{\pm \lambda}(x,-D)$, see also \cite{KN}. Namely, $^{R}\{e^{\pm \lambda}\}(x,D)$ is defined by the oscillatory integral
\begin{eqnarray*}^{R}\{e^{\pm \lambda}\}(x,D)u(x) &=& Os - \iint e^{i\xi (x-y) \pm \lambda(y,\xi)}u(y)\, dy \dslash \xi \\
	&=& \lim_{\varepsilon \to 0} \iint_{\R^{2n}} e^{i\xi (x-y) \pm \lambda(y,\xi)}  e^{-\varepsilon^2y^2} e^{-\varepsilon^2\eta^2} u(y)\, dy \dslash \xi.
\end{eqnarray*}
We observe that when $\lambda$ is a real-valued symbol, the reverse operator $^{R}\{e^{\pm \lambda}\}(x,D)$ coincides with the $L^2$ adjoint $(e^{\pm \lambda}(x,D))^{*}$. The asymptotic expansion of the reverse operator is 
$$
^{R}\{e^{\pm \lambda}\} (x,\xi)\sim \sum_{\alpha} \frac{1}{\alpha!} \partial^{\alpha}_{\xi} D^{\alpha}_{x} e^{\pm\lambda (x,\xi)}.
$$

We have the following continuity result whose proof can be found in Proposition $6.7$ of part I in \cite{KN}. 

\begin{theorem}\label{theorem_continuity_infinite_order_in_Gevrey_sobolev_spaces}
	Let $\rho, m \in \R$ and $\theta, \mu > 1$ with $\theta > 2\mu - 1$. Let $\lambda \in S^{\frac{1}{\theta'}}_{\mu}(\R^{2n})$. Then: \\
	i)
	if $\theta' > \theta$, then $e^{\lambda}(x,D),\, ^{R} (e^{\lambda})(x,D):H^{m}_{\rho; \theta} (\R)\longrightarrow H^{m}_{\rho-\delta; \theta}(\R)$ continuously for every $\delta >0$; \\
	ii) if $\theta' = \theta$, then $e^{\lambda}(x,D),\, ^{R} (e^{\lambda})(x,D):H^{m}_{\rho; \theta} (\R)\longrightarrow H^{m}_{\rho-\delta; \theta}(\R)$ continuously for every 
	$$
	\delta >C(\lambda):= \sup_{(x,\xi) \in \R^{2}} \ds\frac{|\lambda(x,\xi)|}{\langle \xi \rangle^{1/\theta}}.
	$$
\end{theorem}

\section{Proof of Theorem \ref{main}}

We will first prove Theorem $\ref{main}$ in the case $\sigma_2, \sigma_1 \in (0,1)$, since this is the most technical and difficult one. Then, in the final subsections of this section, we shall discuss how to treat the remaining cases: $\sigma_2 \in (0,1)$ and $\sigma_1 \geq 1$, $\sigma_1 \in (0,1)$ and $\sigma_2 \geq 1$, $\sigma_2 \geq 1$ and $\sigma_1 \geq 1$. 

The proof is based on a suitable change of variable of the form 
\beqs\label{ch}
v(t,x)=e^{k(t)\langle D \rangle^{\frac{1}{\theta}}_{h} }e^{\Lambda}(t,x,D)u(t,x)
\eeqs
where $e^{\Lambda}(t,x,D)$ is an invertible pseudodifferential operator of infinite order defined by its symbol $e^{\Lambda(t,x,\xi)}$ and $\langle D \rangle_h$ denotes the Fourier multiplier given by the symbol $\langle\xi\rangle_h:=\sqrt{h^2+\xi^2}$ with the parameter $h \geq 1$.
The order of $\Lambda(t,x,\xi)$ will be $q$ with $0 \leq q < \frac{1}{\theta} < 1$. The leading term of the transformation will be so the operator $e^{k(t)\langle D_x\rangle^{\frac{1}{\theta}}_{h} }$; this means that, looking into the asymptotic expansion of 
$$
P_{k, \Lambda}:= e^{k(t)\langle D \rangle^{\frac{1}{\theta}}_{h}} \circ e^{\Lambda} \circ P\circ (e^{\Lambda})^{-1} \circ e^{-k(t)\langle D \rangle^{\frac{1}{\theta}}_{h}}
$$ 
(that contains the term $ik'(t)\langle D \rangle^{\frac{1}{\theta}}_{h}$ coming out from $e^{k(t)\langle D\rangle_h^{1/\theta}} \circ D_t\circ e^{-k(t)\langle D\rangle_h^{1/\theta}}$) we will be able to absorb errors of order $q$ in the estimates. The function $k(t)$ will be of the form 
\begin{equation}\label{k(t)}
k(t) = \rho_0(T+1-t),
\end{equation}
for some fixed $\rho_0 > 0$.

The main role of this change of variable is to make the transformed Cauchy problem 
\beqs
\label{CPlambda}
\begin{cases}
	P_{k, \Lambda}(t,x,D_t,D_x)v(t,x)=f_{k,\Lambda}(t,x) & (t,x)\in[0,T]\times\R\cr
	v(0,x)=g_{k,\Lambda}(x) & x\in\R
\end{cases}
\eeqs
with 
$$
f_{k,\Lambda}= e^{k(t)\langle D \rangle^{\frac{1}{\theta}}_{h}}e^{\Lambda} f,\quad g_{k,\Lambda}= e^{k(t)\langle D \rangle^{\frac{1}{\theta}}_{h}} e^{\Lambda} g,
$$
well-posed in $L^2(\R)$. 
If we prove well-posedness in $L^2$ for problem \eqref{CPlambda}, then we immediately obtain for \eqref{genCP} the well-posedness in $H^{\infty}_{\theta}(\R)$ for $\theta \in (1, \frac{1}{q})$ by the inversion of formula \eqref{ch}. Indeed, if in \eqref{genCP} we take $f, g \in H^{0}_{\rho;\theta}(\R)$ for some $\rho > 0$, in view of Theorem \ref{theorem_continuity_infinite_order_in_Gevrey_sobolev_spaces} (since $q < \frac{1}{\theta}$) we have
$$
f_{k,\Lambda},g_{k,\Lambda} \in L^{2}(\R)
$$
as long as we choose $\rho_0 < \frac{\rho}{T+1}$. 
Hence there is a unique $v \in L^{2}(\R)$ solution to \eqref{CPlambda}. The inverse of the operator $e^{k(t)\langle D_x \rangle^{\frac{1}{\theta}}_{h}} e^{\Lambda}$ has the following shape 
$$
\{e^{k(t)\langle D_x \rangle^{\frac{1}{\theta}}_{h}} e^{\Lambda}\}^{-1} = (e^{\Lambda})^{-1} e^{-k(t)\langle D_x \rangle^{\frac{1}{\theta}}_{h}}.
$$ 
It turns out (see the upcoming Lemma \ref{lemma_inverse_of_e_Lambda}) that $(e^{\Lambda})^{-1} =\, ^R e^{-\Lambda} \circ r$, for some zero order operator $r$. In this way, again in view of Theorem \ref{theorem_continuity_infinite_order_in_Gevrey_sobolev_spaces}, we have that $u = (e^{k(t)\langle D_x \rangle^{\frac{1}{\theta}}_{h}} e^{\Lambda})^{-1} v$, which is the solution to \eqref{genCP}, belongs to $H^{\infty}_{\theta}(\R)$. Therefore \eqref{genCP} is $H^{\infty}_{\theta}(\R)$ well-posed for all $\theta \in (1,\frac{1}{q})$.

In the following subsections we shall define the symbol $\Lambda(t,x,\xi)$ and then prove that the transformed Cauchy problem \eqref{CPlambda} is well-posed in $L^2(\R)$.

\begin{remark}
	The parameter $h \geq 1$ in the weight $\langle \xi \rangle_{h}$ will be used to obtain several estimates in the paper and also to get the invertibility of the operator $e^{\Lambda}(t,x,D)$. We point out that $\langle \xi \rangle_h$ satisfies
	$$
	|\partial^{\alpha}_{\xi} \langle \xi \rangle_h|\leq C^{\alpha+1}\alpha! \langle \xi \rangle^{-\alpha}_{h},
	$$
	where the constant $C > 0$ does not depend on the parameter $h \geq 1$.
\end{remark}

\subsection{Construction of the symbol $\Lambda$}  

Let us first introduce, for a fixed $\mu>1$ that can be chosen arbitrarily close to $1$, two $\mu-$Gevrey cut-off functions $\psi,\chi \in C^\infty(\R)$ such that
\begin{equation}\label{z}
\psi(y) = \left\lbrace \begin{array}{ll}
	1, & |y|\leq\frac{1}{2} \\ 
	0, & |y|\geq1
\end{array}  \right., \qquad \chi(y) = \left\lbrace \begin{array}{ll}
	1, & |y|\leq\frac{1}{2} \\ 
	0, & |y|\geq1
\end{array}  \right.,\quad y\chi'(y) < 0, \quad  y\in\R,
\end{equation}
$$
|\partial_y^\beta\psi(y)|+|\partial_y^\beta\chi(y)|\leq C^{\beta+1}\beta!^\mu, \quad \beta\in\N_0,
$$
for some positive constant $C$. From now on $C$ will denote a positive constant which may change from line to line. Let moreover $\omega\in C^\infty(\R)$ be a $\mu-$Gevrey function such that 
$$
|\partial_\xi^\alpha\omega(\xi)|\leq C^{\alpha+1}\alpha!^\mu, \quad \alpha\in \N_0,
$$ and 
$$\omega(\xi) = \left\lbrace \begin{array}{ll}
	0, & |\xi|\leq1, \\ 
	-\text{sgn}(a_3(t)), & |\xi|>R,
\end{array}  
\right.$$
for an $R>1$. Given $q_1,q_2\in [0,1)$, we set $q=\max\{q_1,q_2\}$ and define, for some positive constants $M_2,M_1$ to be determined later on, the symbols
\begin{equation}\label{L2}
	\lambda_2(x, \xi) = M_2 w\left(\frac{\xi}{h}\right) \int_{0}^{x} \langle y \rangle^{-\sigma_2}  
	\psi\left(\frac{\langle y \rangle}{\langle \xi \rangle^{\frac{2-q_2}{\sigma_2}}_{h}}\right) dy, \quad (x, \xi) \in \R^2,
\end{equation}
\begin{equation}\label{L1}
	\lambda_1(x, \xi) = M_1 w\left(\frac{\xi}{h}\right) \langle \xi \rangle^{-1}_{h} \int_{0}^{x} \langle y \rangle ^{-\sigma_1} \psi\left(\frac{\langle y \rangle}{\langle \xi \rangle^{\frac{1-q_1}{\sigma_1}}_{h}}\right) dy, \quad (x, \xi) \in \R^2.
\end{equation}
Since on the support of $\psi\left(\langle x \rangle/\langle \xi \rangle^{(2-q_2)/\sigma_2}_{h}\right)$ we have $\langle x \rangle\leq \langle \xi \rangle^{(2-q_2)/\sigma_2}_{h}$ and on the support of $\psi\left(\langle x \rangle/\langle \xi \rangle^{(1-q_1)/\sigma_1}_{h}\right)$ we have $\langle x \rangle\leq \langle \xi \rangle^{(1-q_1)/\sigma_1}_{h}$, the following estimates hold:
\beqs\label{estl2}
|\partial^{\alpha}_{\xi} \partial^{\beta}_{x} \lambda_2(x,\xi)| \leq M_2 \begin{cases}
	C^{\alpha+\beta+1} \alpha!^{\mu} \beta!^{\mu} \langle \xi \rangle^{-\alpha}_{h} \langle x \rangle^{1-\sigma_2-\beta}, \qquad \beta > 0, \alpha \in \N_0,\\
C^{\alpha+1} \alpha!^{\mu}\beta!^{\mu} \langle \xi \rangle^{\frac{(2-q_2)(1-\sigma_2)}{\sigma_2} - \alpha}_{h} \langle x \rangle^{-\beta}, \quad \quad  \beta, \alpha \in \N_0,
	
\end{cases}
\eeqs
\beqs \label{estl1}
|\partial^{\alpha}_{\xi} \partial^{\beta}_{x} \lambda_1(x,\xi)| \leq M_1
	\begin{cases}
		C^{\alpha+\beta+1} \alpha!^{\mu} \beta!^{\mu} \langle \xi \rangle^{-1-\alpha}_{h} \langle x \rangle^{1-\sigma_1-\beta}, \\
		C^{\alpha+\beta+1} \alpha!^{\mu} \beta!^{\mu} \langle \xi \rangle^{\frac{(1-q_1)(1-\sigma_1)}{\sigma_1} - 1 -\alpha}_{h} \langle x \rangle^{-\beta},
	\end{cases}
	\quad \alpha,\beta \in \N_0.
\eeqs

Now we are ready to define the function $\Lambda$ in \eqref{ch}. Taking into account the degeneracy gaps $\ell - k > 0$ and $\ell - k'> 0$, inspired by \cite{CRJEECT}, we define
\begin{equation}\label{defl}
	\Lambda(t,x,\xi) = \Lambda_e(t,x,\xi)+\Lambda_\psi(t,\xi),
\end{equation}
where
\begin{equation}
\Lambda_e(t,x,\xi) = \lambda_{e,2}(t,x,\xi)+\Lambda_{e,2}(t,\xi)+\lambda_{e,1}(t,x,\xi)+\Lambda_{e,1}(t,\xi),
\end{equation}

\beqs\label{lE2}
\lambda_{e,2}(t,x,\xi) &=& t^{-(\ell-k)}\lambda_2(x,\xi) (1-\chi)\left( t\langle\xi\rangle_h^{\frac{2-q_2}{k+1}} \right),
\\
\label{LE2}
\Lambda_{e,2}(t,\xi) &=& -M_{e,2}  \langle \xi \rangle^{(2-q_2)\frac{1-\sigma_2}{\sigma_2}}_{h} \ds \int_0^t \tau^{k-\ell-1} (1-\chi)\left(2\tau\langle\xi\rangle_h^{\frac{2-q_2}{k+1}} \right) d\tau,
\\
\label{lE1}
\lambda_{e,1}(t,x,\xi) &=& t^{-(\ell-k')}\lambda_1(x,\xi) (1-\chi)\left( t\langle\xi\rangle_h^{\frac{1-q_1}{k'+1}} \right),
\\
\label{LE1}
\Lambda_{e,1}(t,\xi) &=& -M_{e,1} \langle \xi \rangle^{\frac{(1-q_1)(1-\sigma_1)}{\sigma_1} - 1}_{h} \ds \int_0^t\tau^{k'-\ell-1} (1-\chi)\left( 2\tau\langle\xi\rangle_h^{\frac{1-q_1}{k'+1}} \right) d\tau,
\eeqs
and
\begin{equation}
	\Lambda_\psi(t,\xi) = \Lambda_{\psi,2}(t,\xi)+\Lambda_{\psi,1}(t,\xi),
\end{equation}
\beqs\label{LP2}
	\Lambda_{\psi,2}(t,\xi) &=& -M_{\psi,2}\langle\xi\rangle_h^2\ds\int_0^t\tau^k\chi \left(\tau\langle\xi\rangle_h^{\frac{2-q_2}{k+1}}\right) d\tau,
\\\
\label{LP1}
\Lambda_{\psi,1}(t,\xi) &=& -M_{\psi,1} \langle\xi\rangle_h\ds\int_0^t\tau^{k'}\chi \left(\tau\langle\xi\rangle_h^{\frac{1-q_1}{k'+1}}\right) d\tau,
\eeqs
with $M_{e,2}, M_{e_1}, M_{\psi,2}, M_{\psi,1}$ positive constants to be chosen later on.
The symbols $\Lambda_e$ and $\Lambda_{\psi}$ are localized, respectively, in the so-called {\it 
evolution zones} $$t\langle\xi\rangle_h^{\frac{2-q_2}{k+1}}\geq 1/4,\quad  t\langle\xi\rangle_h^{\frac{1-q_1}{k'+1}}\geq 1/4$$ and {\it pseudodifferential zones} $$ t\langle\xi\rangle_h^{\frac{2-q_2}{k+1}}\leq1 \quad  t\langle\xi\rangle_h^{\frac{1-q_1}{k'+1}}\leq 1.$$

In view of the support of the cutoff functions $\chi$ and $1-\chi$ and using \eqref{estl2}, \eqref{estl1}, we have the following estimates (where all the constants are independent of time and of the parameter $h$):
\beqs\label{g2}
|\partial^{\alpha}_{\xi} \partial^{\beta}_{x} \lambda_{e,2}(t,x,\xi)|& \leq &M_2 \begin{cases}
	C^{\alpha+1} \alpha!^{\mu} \langle \xi \rangle^{(2-q_2)\left( \frac{\ell-k}{k+1} + \frac{1-\sigma_2}{\sigma_2} \right)-\alpha}_{h},\qquad \quad \quad \beta = 0, \\
	C^{\alpha+\beta+1} \alpha!^{\mu} \beta!^{\mu} \langle \xi \rangle^{(2-q_2)\frac{\ell-k}{k+1}-\alpha}_{h} \langle x \rangle^{1-\sigma_2-\beta}, \,\,\, \beta > 0,
\end{cases}
\\
|\partial^{\alpha}_{\xi} \Lambda_{e,2}(t,\xi)| &\leq &M_{e,2} C^{\alpha+1} \alpha!^{\mu} \langle \xi \rangle^{(2-q_2)\left( \frac{\ell-k}{k+1} + \frac{1-\sigma_2}{\sigma_2} \right)-\alpha}_{h},
\\
|\partial^{\alpha}_{\xi} \Lambda_{\psi,2}(t,\xi)| &\leq & M_{\psi,2} C^{\alpha+1} \alpha!^{\mu} \langle \xi \rangle^{q_2-\alpha}_{h},
\eeqs	
\beqs\label{g1}
|\partial^{\alpha}_{\xi} \partial^{\beta}_{x} \lambda_{e,1}(t,x,\xi)| &\leq& M_1
\begin{cases}
	C^{\alpha+\beta+1} \alpha!^{\mu} \beta!^{\mu} \langle \xi \rangle^{(1-q_1) \frac{\ell-k'}{k'+1}-1-\alpha}_{h}\langle x\rangle^{1-\sigma_1-\beta}, 
	\\
	C^{\alpha+\beta+1} \alpha!^{\mu} \beta!^{\mu} \langle \xi \rangle^{(1-q_1)\left( \frac{1-\sigma_1}{\sigma_1} + \frac{\ell-k'}{k'+1} \right)-1-\alpha}_{h} \langle x \rangle^{-\beta},
\end{cases}
\\\nonumber
|\partial^{\alpha}_{\xi} \Lambda_{e,1}(t,\xi)| &\leq &M_{e,1} C^{\alpha+1} \alpha!^{\mu} \langle \xi \rangle^{(1-q_1)\left( \frac{1-\sigma_1}{\sigma_1} + \frac{\ell-k'}{k'+1} \right)-1-\alpha}_{h},
\\
|\partial^{\alpha}_{\xi} \Lambda_{\psi,1}(t,\xi)| &\leq& M_{\psi,1} C^{\alpha+1} \alpha!^{\mu} \langle \xi \rangle^{q_1-\alpha}_{h}.
\eeqs
Since $\Lambda$ is given by the sum of $\Lambda_e$ and $\Lambda_\psi$, the order of $\Lambda_e$ is a decreasing function of $q_2$ and $q_1$ and the order of $\Lambda_\psi$ is an increasing function of $q_2$ and $q_1$, to minimize the total order of $\Lambda$ we take  $q_1$ and $q_2$ satisfying
\beqs\label{choiceq12}
\left\{
\begin{array}{l}
q_2 = (2-q_2)\left( \frac{\ell-k}{k+1} + \frac{1-\sigma_2}{\sigma_2} \right),
\\
q_1 = (1-q_1)\left( \frac{\ell-k'}{k'+1} + \frac{1-\sigma_1}{\sigma_1} \right)-1.
\end{array}
\right.
\eeqs
Hence, we have that the pseudodifferential and evolution zones have the same influence on the order of the symbol $\Lambda$. Solving the above two equations we obtain
\beqs\label{q2}
q_2 = 2\frac{(\ell-k)\sigma_2 + (k+1)(1-\sigma_2)}{\sigma_2(\ell-k)+(k+1)}
= 2 \left( 1 - \frac{\sigma_2(k+1)}{\sigma_2(\ell-k)+(k+1)}\right),
\eeqs
\beqs\label{q1}
q_1 = \frac{(k'+1)(1-2\sigma_1)+\sigma_1(\ell-k')}{\sigma_1(\ell-k')+(k'+1)}
= 1 - \frac{2\sigma_1(k'+1)}{\sigma_1(\ell-k') + (k'+1)}.
\eeqs
Observe that on the one hand $q_1 < 1$, while, if $\ell < 2+3k$ (which is true since we are working under the condition $\ell<2k+1$), then
$$
q_2 < 1 \iff \sigma_2 > \frac{k+1}{2+3k-\ell} = \frac{k+1}{(k+1) + (1+2k-\ell)},
$$
that is assumption \eqref{lowbound}.
We also point out that $q_2 > 0$ for all $\sigma_2 \in (0,1)$ and $\ell > k$. The total order of $\Lambda$ will be then $q=\max\{q_1,q_2\}\in (0,1)$.


Whith these choices of $q_1,q_2$ and using \eqref{choiceq12} we get the following:
\beqs\nonumber
\Lambda_2(t,x,\xi):&=&\lambda_{e,2}(t,x,\xi)+\Lambda_{e,2}(t,\xi)+\Lambda_{\psi,2}(t,\xi)\in SG_\mu^{q_2,0}(\R^2),
\\
\nonumber
\Lambda_1(t,x,\xi):&=&\lambda_{e,1}(t,x,\xi)+\Lambda_{e,1}(t,\xi)+\Lambda_{\psi,1}(t,\xi)\in SG_\mu^{q_1,0}(\R^2),
\\
\nonumber
\Lambda(t,x,\xi)&=& \Lambda_2(t,x,\xi)+\Lambda_1 (t,x,\xi)\in SG^{q,0}_{\mu}(\R^2),
\eeqs
and moreover using \eqref{g2}, \eqref{g1}, \eqref{q2}, and \eqref{q1}, we have
\beqs
\label{d1}
&&\partial_x\Lambda_2=\partial_{x} \lambda_{e,2} \in SG^{(2-q_2)\frac{\ell-k}{k+1}, -\sigma_2}_{\mu}(\R^2) =  
SG^{q_2 - \frac{2(k+1)(1-\sigma_2)}{\sigma_2(\ell-k)+(k+1)}, -\sigma_2}_{\mu}(\R^2),
\\\label{d2}
&&\partial_x\Lambda_1=\partial_x \lambda_{e,1}\in  SG^{(1-q_1)\frac{\ell-k'}{k'+1} - 1, -\sigma_1}_{\mu}(\R^2) = 
SG^{q_1 - \frac{2(k'+1)(1-\sigma_1)}{\sigma_1(\ell-k')+(k'+1)}, -\sigma_1}_{\mu}(\R^2).
\eeqs

We conclude this section with some remarks on the structure of the time derivate of $\Lambda(t,x,\xi)$:
{\small \beqsn
		\partial_t&& \hskip-0.5cm\Lambda (t,x,\xi)  
\\
&=& t^{-(\ell-k) - 1} \left\{ -(\ell-k) \lambda_2 (1-\chi)\left( t\langle\xi\rangle_h^{\frac{2-q_2}{k+1}} \right) 
		-\lambda_2 \chi'\left( t\langle\xi\rangle_h^{\frac{2-q_2}{k+1}} \right)t\langle\xi\rangle_h^{\frac{2-q_2}{k+1}}
		- M_{e,2} \langle \xi \rangle^{(2-q_2)\frac{1-\sigma_2}{\sigma_2}}_{h}  \right\}
\\
&&\cdot(1-\chi)\left( 2t\langle\xi\rangle_h^{\frac{2-q_2}{k+1}} \right)- M_{\psi,2}\langle\xi\rangle_{h}^{2} t^k \chi \left(t\langle\xi\rangle_h^{\frac{2-q_2}{k+1}}\right)
 \\
		&+& t^{-(\ell-k') - 1} \left\{ -(\ell-k') \lambda_1 (1-\chi)\left( t\langle\xi\rangle_h^{\frac{1-q_1}{k'+1}} \right) 
		-\lambda_1 \chi'\left( t\langle\xi\rangle_h^{\frac{1-q_1}{k'+1}} \right)t\langle\xi\rangle_h^{\frac{1-q_1}{k'+1}}
		- M_{e,1} \langle \xi \rangle^{(1-q_1)\frac{1-\sigma_1}{\sigma_1}}_{h}  \right\}
\\
&&\cdot(1-\chi)\left( 2t\langle\xi\rangle_h^{\frac{1-q_1}{k'+1}} \right) - 
		 M_{\psi,1} \langle\xi\rangle_h t^{k'} \chi \left(t\langle\xi\rangle_h^{\frac{1-q_1}{k'+1}}\right)
\\
&=& a(\Lambda_2)(t,x,\xi) + a(\Lambda_1)(t,x,\xi),
\eeqsn}
where on the one hand we have
\beqs\label{star2}
|a(\Lambda_2)(t,x,\xi)| &\lesssim& \langle \xi \rangle^{q_2 + \frac{2-q_2}{k+1}}_{h} = \langle \xi \rangle^{2-\frac{(2-q_2)k}{k+1}}_{h},
\\ \label{star1}
|a(\Lambda_1)(t,x,\xi|)& \lesssim& \langle \xi \rangle^{q_1 + \frac{1-q_1}{k'+1}}_{h} = \langle \xi \rangle^{1-\frac{(1-q_1)k'}{k'+1}}_{h},
\eeqs
with uniform bounds with respect to time, that is, $a(\Lambda_2)$ has order strictly less than $2$ and $a(\Lambda_1)$ has order strictly less than $1$, and on the other hand
\beqs\label{dta2}
|a(\Lambda_2)(t,x,\xi)| &\lesssim &t^{k} \langle \xi \rangle^{2}_{h},
\\\label{dta1}
| a(\Lambda_1)(t,x,\xi)| &\lesssim &t^{k'} \langle \xi \rangle^{1}_{h}.
\eeqs


\subsection{Invertibility of the operator $e^{\Lambda}$}

Since $\Lambda(t,x,\xi) \in SG^{q,0}_{\mu}(\R^2)$, by Proposition \ref{proposition_exponential_of_symbols_of_finite_order} we have $e^{\Lambda(t,x,\xi)} \in SG^{\infty, 0}_{\mu; \frac{1}{q}}(\R^2)$. The calculus for this class of symbols (see Section \ref{preliminaries})  is defined for $1/q>1$. This is the reason for the restriction \eqref{lowbound} in Theorem \ref{main}, see also Remark \ref{rem1}.

Since $\Lambda(t,x,\xi)$ is real-valued the reverse operator $^Re^{-\Lambda}$ coincides with the $L^2$-adjoint of $e^{-\Lambda}$ and by Theorem \ref{theorem_symbolic_calculus_of_infinte_order} we have
\beqsn
&&^{R}e^{-{\Lambda}} (t,x,D)= a(t,x,D) + r(t,x,D),
\\
 &&a (t,x,\xi)\sim \sum_{\alpha} \frac{1}{\alpha!} \partial^{\alpha}_{\xi}D^{\alpha}_{x}e^{-{\Lambda(t,x,\xi)}}\ {\rm in}\ FSG^{\infty, 0}_{\mu; 1/q}(\R^{2}),\; r \in \mathcal{S}_{2\mu - 1}(\R^{2}).
\eeqsn
Therefore (again by Theorem \ref{theorem_symbolic_calculus_of_infinte_order}) we get 
\beqsn
	(e^\Lambda \circ ^{R}e^{-\Lambda}) (t,x,D)= q(t,x,D) + r_{\infty}(t,x,D),
\eeqsn
where $r_{\infty}(t,x,\xi) \in \mathcal{S}_{2\mu-1}$ and 
$$
q(t,x,\xi) \sim \sum_{\alpha \geq 0}  \frac{1}{\alpha!}  \partial^{\alpha}_{\xi} \{e^{\Lambda(t,x,\xi)}  D^{\alpha}_{x}e^{-\Lambda(t,x,\xi)} \} \quad \text{in} \quad FSG^{\infty,0}_{\mu;\frac{1}{q}}(\R^2).
$$
Standard computations give that 
$$
|\partial^{\gamma+\alpha}_{\xi} \partial^{\beta}_{x} \{e^{\Lambda(t,x,\xi)}  D^{\alpha}_{x}e^{-\Lambda(t,x,\xi)} \}| \leq 
C^{\gamma+\beta+2\alpha+1} (\gamma!\beta!\alpha!^{2})^{\mu} \langle \xi \rangle^{-(1-q)\alpha - \gamma}_{h} \langle x \rangle^{-\beta-\alpha}.
$$
In this way we obtain that $q (t,x,\xi)\in SG^{0,0}_{\mu}(\R^2)$ and 
$$
 (e^\Lambda \circ ^{R}e^{-\Lambda}) (t,x,D)= q(t,x,D) + r_{\infty}(t,x,D) = 1 + r(t,x,D) + r_{\infty}(t,x,D),
$$
where $r (t,x,\xi)\in SG^{-(1-q), -1}_{\mu}(\R^2)$ and $r_{\infty}(t,x,\xi) \in \mathcal{S}_{2\mu-1}(\R^2)$. Setting $\bar{r} = r+r_{\infty}$ we get 
\beqs\label{v}
(e^\Lambda \circ ^{R}e^{-\Lambda}) (t,x,D) = 1 + \bar{r}(t,x,D),
\eeqs
and the $0-$seminorms of the symbol $\bar{r}$ can be estimated by 
\begin{equation}\label{silksong}
|\partial^{\alpha}_{\xi} \partial^{\beta}_{x} \bar{r}(t,x,\xi)| \leq C_{\alpha,\beta} \langle \xi \rangle^{-(1-q) - \alpha}_{h} \leq h^{-(1-q)} C_{\alpha,\beta} \langle \xi \rangle^{-\alpha}_{h}.
\end{equation}
Calderon-Vaillancourt theorem then implies that the $L^2(\R)$ operator norm of $\bar{r}$ is strictly less than $1$, provided that $h \geq h_0$ for some $h_0$ depending on a finite number of the constants $C_{\alpha, \beta}$ appearing in \eqref{silksong}. Therefore, the operator $1 + \bar{r}$ is invertible on $L^2(\R)$ and we have 
$$
(1+\bar{r}) \circ \sum_{j = 0}^{\infty} (-\bar{r})^{j} = 1.
$$ 
Notice that since the constants $C_{\alpha,\beta}$ depend on the constants $M_2$ and $M_1$, then $h_0$ depends on $M_2$ and $M_1$ as well. Finally we get by \eqref{v}
$$
e^\Lambda \circ ^{R}e^{-\Lambda} \circ \sum_{j = 0}^{\infty} (-\bar{r})^{j} = 1.
$$

Using Theorem $9$ in \cite{AAC3evolGelfand-Shilov} we obtain that $\sum_{j = 0}^{\infty} (-\bar{r})^{j}$ is an operator given by a symbol of the form $q_1 + \bar r_{\infty}$, where $q_1(t,x,\xi) \in SG^{0,0}_{\mu}$ and $\bar{r}_{\infty} (t,x,\xi)\in \Sigma_{\kappa}$ for every $\kappa > 2\mu -1$. We have so proved the following lemma:

\begin{lemma}\label{lemma_inverse_of_e_Lambda}
	Let $\mu > 1$. There is $h_0 = h_0(M_2,M_1)$ such that for all $h \geq h_0$, the operator $e^{\Lambda}(t,x,D)$ is invertible 
	and its inverse is given by 
	$$
	\{e^{\Lambda}(t,x,D)\}^{-1} = \hskip2pt ^R  \{e^{-\Lambda}\}(t,x,D) \circ \sum_{j \geq 0} (-\overline{r}(t,x,D))^{j}, 
	$$
	for some $\bar{r} = r + r_{\infty}$, where $r (t,x,\xi)\in SG^{-(1-q),-1}_{\mu}(\R^{2})$, $r_{\infty} (t,x,\xi)\in \Sigma_{\kappa}(\R^{2})$ for every 
	$\kappa> 2\mu-1$ and 
	$$
	r (t,x,\xi)- \sum_{1 \leq \gamma \leq N} \frac{1}{\gamma!} \partial^{\gamma}_{\xi} \{ e^{\Lambda(t,x,\xi)} D^{\gamma}_{x} e^{-\Lambda(t,x,\xi)} \} \in SG^{-(N+1)(1-q),-N-1}_{\mu}(\R^{2}), \quad \forall N \geq 1.
	$$ 
	Moreover, $\sum (-\overline{r}(t,x,D))^{j}$ has symbol in $SG^{0,0}_{\mu}(\R^{2}) + \Sigma_{\kappa}(\R^{2})$, for every $\kappa >2\mu-1$. Finally, we have
		\begin{equation}\label{equation_inverse_of_e_power_tilde_Lambda_in_a_precise_way}
		\{e^{\Lambda}(t,x,D) \}^{-1} = \hskip2pt ^R \{e^{-\Lambda}\}(t,x,D) \circ \textrm{op} ( 1 - i\partial_{\xi} \partial_{x} {\Lambda} + q_{-2} + r_{\infty} ),
	\end{equation}
	where $q_{-2}(t,x,\xi) \in SG^{-2(1-q),-2}_{\mu}(\R^{2})$ and  $r_{\infty}(t,x,\xi) \in \Sigma_{\kappa}(\R^{2})$, for all $\kappa > 2\mu-1$.
\end{lemma}

\subsection{Conjugation $P_{\Lambda}:=(e^\Lambda)P(e^{\Lambda})^{-1}$} 

This section is devoted to compute the conjugation$$P_{\Lambda}=(e^\Lambda)\circ P\circ (e^{\Lambda})^{-1}.$$ The coefficients of the operator $P$ belong to the class $S^{j}_{1}(\R^2)$, for $j =0,1,2,3$, so we start by computing $(e^\Lambda) p (e^{\Lambda})^{-1}$ where  $p(t,x,\xi) \in S^{m}_{\mu}(\R^2)$, $m \in \R$. Since the inverse of $e^{\Lambda}$ contains the reverse operator we first calculate
$$
e^{\Lambda} (t,x,D)\circ p(t,x,D) \circ\, ^{R} e^{-\Lambda} (t,x,D)= q(t,x,D) + r_{\infty}(t,x,D),
$$
where $r_\infty (t,x,\xi)\in \mathcal{K}_{2\mu-1}$ and 
$$
q(t,x,\xi) \sim \sum_{\alpha,\beta} \frac{1}{\alpha!\beta!} \partial^{\alpha}_{\xi} \{\partial^{\beta}_{\xi} e^{\Lambda(t,x,\xi)} D^{\beta}_{x}p(t,x, \xi) D^{\alpha}_{x}e^{-\Lambda(t,x,\xi)} \} \quad \text{in} \quad FS^{\infty}_{\mu}(\R^2).
$$
Standard computations at symbols level give 
$$
|\partial^{\alpha+\gamma}_{\xi} \partial^{\lambda}_{x} \{\partial^{\beta}_{\xi} e^{\Lambda} D^{\beta}_{x}p D^{\alpha}_{x}e^{-\Lambda}\}| \leq
C^{\gamma+\lambda+2(\alpha+\beta)+1} (\gamma!\lambda!\alpha!^{2}\beta!^{2})^{\mu} \langle \xi \rangle^{m-(1-q)(\alpha+\beta) - \gamma}_{h} .
$$
Therefore $q(t,x,\xi) \in S^{m}_{\mu}(\R^2)$
and we can write at operators level
$$
e^{\Lambda} \circ p \circ \, ^{R} e^{-\Lambda} = \sum_{0 \leq \alpha + \beta < N} \frac{1}{\alpha!\beta!} op\{\partial^{\alpha}_{\xi} \{\partial^{\beta}_{\xi} e^{\Lambda} D^{\beta}_{x}p D^{\alpha}_{x}e^{-\Lambda} \} \} + r_{N} + r_{\infty},
$$
where $r_{N}(t,x,\xi) \in S^{m-N(1-q)}_{\mu}(\R^2)$ and $r_{\infty} (t,x,\xi)\in \mathcal{K}_{2\mu-1}$. In particular, if $p$ satisfies 
$$
|\partial^{\alpha}_{\xi} \partial^{\beta}_{x} p(t,x,\xi)| \leq t^{s}C^{\alpha+\beta+1} \alpha!^{\mu}\beta!^{\mu} \langle \xi \rangle^{m-\alpha}_{h} \langle x \rangle^{-\sigma},
$$
then 
\begin{equation}
(e^{\Lambda} \circ p \circ \, ^{R} e^{-\Lambda})(t,x,D) = p(t,x,D) + p_{q-1}(t,x,D) + r_{\infty}(t,x,D),
\end{equation}
where $r_{\infty}(t,x,\xi) \in \mathcal{K}_{2\mu-1}$ and $p_{q-1} (t,x,\xi)\in S^{m-(1-q)}_{\mu}(\R^2)$ satisfies
\begin{equation}\label{lace}
|\partial^{\alpha}_{\xi} \partial^{\beta}_{x} p_{q-1}(t,x,\xi)| \leq t^s C^{\alpha+\beta+1} \alpha!^{\mu}\beta!^{\mu} \langle \xi \rangle^{m-(1-q)-\alpha}_{h} \langle x \rangle^{-\sigma}.
\end{equation}
Finally, composing with the Neumann series we get (for new operators $p_{q-1}$ and $r_{\infty}$)
\begin{equation}\label{hornet}
	(e^{\Lambda} \circ p \circ (e^{\Lambda})^{-1})(t,x,D) = p (t,x,D)+ p_{q-1} (t,x,D)+ r_{\infty}(t,x,D),
\end{equation}
where  $r_{\infty}(t,x,\xi) \in \Sigma_{\kappa}(\R^2)$ for every $\kappa > 2\mu-1$, $p_{q-1}(t,x,\xi) \in S^{m-(1-q)}_{\mu}(\R^2)$ and $p_{q-1}(t,x,\xi)$ still satisfies \eqref{lace}. 

\bigskip

\noindent{\textbf{- Conjugation of} $D_t$}:
Note that 
$$
D_t \circ e^{\Lambda(t)} = e^{\Lambda(t)} \circ D_t + op\{ e^{\Lambda(t)}D_t\Lambda\}
$$
hence
\beqsn
e^{\Lambda(t)} \circ D_t \circ \{e^{\Lambda(t)}\}^{-1} &=& D_t - op\{ e^{\Lambda(t)}D_t\Lambda  \}  \{e^{\Lambda(t)}\}^{-1}
\\
&=&D_t - op\{ e^{\Lambda(t)}D_t\Lambda  \}\circ ^R \{e^{-\Lambda}(t,x,D)\} op( 1 - i\partial_{\xi} \partial_{x} {\Lambda} + q_{-2} + r_{\infty} )
\eeqsn
by \eqref{equation_inverse_of_e_power_tilde_Lambda_in_a_precise_way}. We have 
\begin{align*}
	op\{ e^{\Lambda(t)}D_t\Lambda\} &\circ \,^R\{e^{-\Lambda(t)}\} =  q_1 + r_{\infty},
\end{align*}
where $r_{\infty} \in \mathcal{S}_{2\mu-1}(\R^2)$ and 
$$
q_1 \sim \sum_{\gamma} \frac{1}{\gamma!} \partial^{\gamma}_{\xi} \{e^{\Lambda(t)} D_t\Lambda(t)D^{\gamma}_{x} e^{-\Lambda(t)}\} \quad \text{in} \quad FSG^{\infty,0}_{\mu}(\R^2).
$$
In fact, since $\partial_t \Lambda \in SG^{2,0}_{\mu}(\R^2)$ we have
$$
|\partial^{\gamma+\alpha}_{\xi} \partial^{\beta}_{x} \{e^{\Lambda(t)} D_t\Lambda(t)D^{\gamma}_{x} e^{-\Lambda(t)}\}| \leq C^{\alpha+\beta+2\gamma+1} (\alpha!\beta!\gamma!^{2})^{\mu} \langle \xi \rangle^{2-\alpha - \gamma(1-q)}_{h} \langle x \rangle^{-\beta-\gamma}.
$$
Therefore $q_1 \in SG^{2,0}_{\mu}(\R^2)$.
Writing $\partial_{t} \Lambda = a(\Lambda_2) + a(\Lambda_1)$ and recalling that $a(\Lambda_2)$ has order $2$, $a(\Lambda_1)$ has order $1$, and $\Lambda\in SG^{q,0}_{\mu}$ we get 
$$
q_1 (t,x,D)= D_t \Lambda (t,x,D)+ r_{2,D_t} (t,x,D)+ r_{1,D_t}(t,x,D), 
$$
where $r_{2,D_t}(t,x,\xi) \in SG^{2-(1-q), -1}_{\mu} (\R^2)$ and $r_{1,D_t} (t,x,\xi)\in SG^{1-(1-q), -1}_{\mu} (\R^2)$. Moreover, taking \eqref{dta2} and \eqref{dta1} into account we obtain
$$
|\partial^{\alpha}_{\xi} \partial^{\beta}_{x} r_{2,D_t}(t,x,\xi)| \leq t^{k} C^{\alpha+\beta+1}(\alpha!\beta!)^{\mu} \langle \xi \rangle^{2-(1-q)}_{h} \langle x \rangle^{-1},$$
$$
|\partial^{\alpha}_{\xi} \partial^{\beta}_{x} r_{1,D_t}(t,x,\xi)| \leq t^{k'} C^{\alpha+\beta+1}(\alpha!\beta!)^{\mu} \langle \xi \rangle^{1-(1-q)}_{h} \langle x \rangle^{-1}.
$$
Composing with the Neumann series the same structure remains, more precisely we have
\begin{align*}
	\{&D_t\Lambda e^{\Lambda}\} \{e^{\Lambda}\}^{-1}(t,x,D) \\
	&= \{D_t\Lambda (t,x,D)+ r_{2,D_t} (t,x,D)+ r_{1,D_t}(t,x,D) \}op \{ 1 - i\partial_\xi\partial_x \Lambda + q_{-2} + r_{\infty} \} \\
	&= D_{t}\Lambda(t,x,D) + r_{D_t} (t,x,D)+ \tilde{r}_{D_t}(t,x,D) + r_{\infty}(t,x,D), 
\end{align*}
where $r_{D_t} (t,x,\xi)\in SG^{2 - (1-q), -1 }_{\mu}$, $\tilde{r}_{D_t} (t,x,\xi)\in SG^{1-(1-q), -1}_{\mu},$ $r_{\infty}(t,x,\xi) \in \Sigma_{\kappa}(\R^2)$ for all $\kappa > 2\mu-1$ and 
\begin{equation}\label{melanie_1}
|\partial^{\alpha}_{\xi} \partial^{\beta}_{x} r_{D_t}(t,x,\xi)| \leq t^{k} C^{\alpha+\beta+1} (\alpha!\beta!)^{\mu} \langle \xi \rangle^{2-(1-q)}_{h} \langle x \rangle^{-1}, 
\end{equation}
\begin{equation}\label{melanie_2}
|\partial^{\alpha}_{\xi} \partial^{\beta}_{x} \tilde{r}_{D_t}(t,x,\xi)| \leq t^{k'} C^{\alpha+\beta+1} (\alpha!\beta!)^{\mu} \langle \xi \rangle^{1-(1-q)}_{h} \langle x \rangle^{-1}.	
\end{equation}

\bigskip

\noindent\textbf{- Conjugation of} $a_3(t)D^3_{x}$.
We have
\begin{align*}
	e^\Lambda (t,x,D) \circ a_3(t)D^{3}_{x} \circ {^R}e^{-\Lambda}  (t,x,D) &= q_1(t,x,D )+ r_{\infty}(t,x,D),
\end{align*}
where $r_{\infty}(t,x,\xi)\in \mathcal{S}_{2\mu-1}(\R)$ and
$$
q_1 (t,x,\xi)\sim  \sum_{\gamma} \frac{1}{\gamma!} \partial^{\gamma}_{\xi} \{e^{\Lambda(t,x,\xi)} a_3(t) \xi^3 D^{\gamma}_{x}e^{-\Lambda(t,x,\xi)} \} \quad \text{in}  \quad FSG^{3,0}_{\mu}(\R^2).
$$
Let us observe that since we are assuming $\eqref{a3}$, that is $|a_3(t)|$ is comparable with $t^{\ell}$, when we multiply $\partial_x\Lambda$ by $a_3(t)$ we get rid of the negative powers of $t$ appearing in $\partial_x \Lambda$, hence
\beqs\label{peach}
a_3(t)\partial_{x}\Lambda (t,x,\xi)&\sim &t^{k} \partial_{x}\lambda_2 (1-\chi)\left( t\langle\xi\rangle_h^{\frac{2-q_2}{\ell+1}} \right) 
+ t^{k'} \partial_{x}\lambda_1 (1-\chi)\left( t\langle\xi\rangle_h^{\frac{1-q_1}{\ell+1}} \right)
\\\nonumber
&\in &SG^{0, -\sigma_2}_{\mu}(\R^2) + SG^{-1,-\sigma_1}_{\mu}(\R^2).
\eeqs
Thus, using \eqref{peach} we get
$$
 (e^\Lambda   \circ a_3(t) D^{3}_{x} \circ {^R}e^{-\Lambda} ) (t,x,D) = op(a_3(t) \xi^3 + i a_3(t)3\xi^2\partial_x \Lambda + ia_3(t)\xi^3 \partial_\xi \partial_x \Lambda + d_1 + r_0 + r_{\infty})
$$
where $d_1(t,x,\xi)$ is bounded by $t^{k}$, does not depend on $\lambda_1$ and belongs to $SG^{1,-\sigma_2-1}_{\mu}(\R^2)$ and $r_0$ is of order zero.

Now we have to compose with the Neumann series. We have that 
\beqsn
e^\Lambda&&\hskip-0.2cm(t,x,D(t,x,D)) \circ a_3(t) D^{3}_{x} \circ (e^{\Lambda})^{-1} = a_3(t)D_x^3\circ \sum_{j = 0}^{\infty}  (-\overline{r})^{j} (t,x,D)
\\
&&+op(i a_3(t)3\xi^2\partial_x \Lambda + ia_3(t)\xi^3 \partial_\xi \partial_x \Lambda + d_1 + r_0 + r_{\infty})\circ op \{ 1 - i\partial_\xi\partial_x \Lambda + q_{-2} + r_{\infty} \}.
\eeqsn
Concerning the first term on the right hand side we have (using \eqref{peach})
$$
a_3(t) \sum_{j = 0}^{\infty}  (-\overline{r})^{j} (t,x,D)= a_3(t)op(1-i\partial_{\xi}\partial_{x}\Lambda+ q_{-1} + r_{-3} + r_{\infty}),
$$
where $a_3(t)q_{-1}(t,x,\xi) \in SG^{-2,-\sigma_2-1}_{\mu}(\R^2)$, $a_3(t)q_{-1}$ does not depend on $\lambda_1$, $a_3(t)r_{-3}$ denotes an operator of order -3 with bounded in time coefficients and $a_3(t)r_{\infty} \in \Sigma_{\kappa}(\R^2)$ for all $\kappa > 2\mu-1$. Therefore 
$$
a_3(t) D_x^3 \sum_{j \geq 0} (-\overline{r})^{j} (t,x,D)= a_3(t) D_x^3 - ia_3(t) op(\xi^3 \partial_{\xi}\partial_{x} \Lambda) + d_1(t,x,D) + r_0(t,x,D) + r_{\infty}(t,x,D),
$$
for some new operators $d_1$, and $r_{\infty}$ satisfying the same properties as before, and $r_0$ of order zero with bounded in time coefficients. For the second term, note that the symbol of the composed operator $$op(i a_3(t)3\xi^2\partial_x \Lambda + ia_3(t)\xi^3 \partial_\xi \partial_x \Lambda + d_1 + r_0 + r_{\infty})\circ op\{ 1 - i\partial_\xi\partial_x \Lambda + q_{-2} + r_{\infty} \}$$ is given by
\begin{align*}
	 ia_3(t)3\xi^2\partial_x \Lambda + a_{2-(1-q)} + ia_3(t)op(\xi^3 \partial_\xi \partial_x \Lambda)  + a_{1-(1-q)} + r_0 + r_{\infty},
\end{align*}
where $a_{2-(1-q)} \in SG^{2-(1-q), -1-\sigma_2}_{\mu}$ contains also the term $d_1$ and is bounded by $t^{k}$, where $a_{1-(1-q)} \in SG^{1-(1-q), -1-\sigma_1}_{\mu}$  is bounded by $t^{k'}$, and  where $r_0$, $r_{\infty}$ are as before. Finally, we get 
\begin{align*}
	(e^{\Lambda} a_3(t) D^{3}_{x} \{e^{\Lambda}\}^{-1} )(t,x,D)
	&= a_3(t) D_x^3 + ia_3(t)op(3\xi^2\partial_x \Lambda) + (r_{a_3} + \tilde{r}_{a_3} + r_0+ r_{\infty})(t,x,D),
\end{align*}
where $r_{a_3}(t,x,\xi) \in SG^{2 - (1-q), -1-\sigma_2}_{\mu}$, $\tilde{r}_{a_3}(t,x,\xi) \in SG^{1-(1-q), -1-\sigma_1}_{\mu}$ and
$$
|\partial^{\alpha}_{\xi} \partial^{\beta}_{x} r_{a_3}(t,x,\xi)| \leq t^{k} C^{\alpha+\beta+1} (\alpha!\beta!)^{\mu} \langle \xi \rangle^{2-(1-q)}_{h} \langle x \rangle^{-1-\sigma_2}, $$
$$|\tilde{r}_{a_3}(t,x,\xi)| \leq t^{k'} C^{\alpha+\beta+1} (\alpha!\beta!)^{\mu} \langle \xi \rangle^{1-(1-q)}_{h} \langle x \rangle^{-1-\sigma_1}.
$$

\bigskip

\noindent\textbf{- Conjugation of} $a_2(t,x)D^{2}_{x}$.
We just apply formula \eqref{hornet} to get
$$
e^\Lambda(t,x,D) a_2(t,x)D^{2}_{x}\, \{e^{\Lambda}\}^{-1}(t,x,D) = a_2(t,x)D^{2}_x +( r_{a_2} + r_{\infty})(t,x,D),
$$
where $r_{a_2}$ satisfies
\begin{equation}\label{melanie_3}
|\partial^{\alpha}_{\xi} \partial^{\beta}_{x} r_{a_2}(t,x,\xi)| \leq t^{k}C^{\alpha+\beta+1} \alpha!^{\mu}\beta!^{\mu} \langle \xi \rangle^{2-(1-q)-\alpha}_{h} \langle x \rangle^{-\sigma_2}.
\end{equation}

\bigskip

\noindent \textbf{- Conjugation of} $a_1(t,x)D_{x}$.
Again we just apply formula \eqref{hornet} to  get
$$
e^\Lambda(t,x,D) a_1(t,x) D_x \, \{e^{\Lambda}\}^{-1} (t,x,D)= a_1(t,x)D_x + (r_{a_1} + r_0 + r_{\infty})t,x,D),
$$
where $r_{a_1}$ satisfies
\begin{equation}\label{melanie_4}
|\partial^{\alpha}_{\xi} \partial^{\beta}_{x} r_{a_1}(t,x,\xi)| \leq t^{k'}C^{\alpha+\beta+1} \alpha!^{\mu}\beta!^{\mu} \langle \xi \rangle^{1-(1-q)-\alpha}_{h} \langle x \rangle^{-\sigma_1}.
\end{equation}

\bigskip

Summing up we get the following expression for the operator $P_{\Lambda}$:
\begin{align*}
	P_{\Lambda}(t,x,D) = (e^{\Lambda} P \{e^{\Lambda}\}^{-1})(t,x,D) &= D_t - D_t\Lambda (t,x,D) \\ 
	&+ a_3(t) D^3_x + ia_3(t) op(3\xi^2 \partial_x \Lambda) + r_{a_3} (t,x,D)+ \tilde{r}_{a_3}(t,x,D) \\
	&+ a_2(t,x)D^{2}_{x} + r_{a_2}(t,x,D) \\
	&+ a_1(t,x)D_{x} + r_{a_1}(t,x,D) + r_{0}(t,x,D) + r_{\infty}(t,x,D),
\end{align*}
where we have absorbed the terms $r_{D_t}$ and $\tilde{r}_{D_t}$ into $r_{a_2}$ and $r_{a_1}$, respectively (cf. \eqref{melanie_1}, \eqref{melanie_2}, \eqref{melanie_3} and \eqref{melanie_4}).  

\subsection{Conjugation $P_{k, \Lambda} = e^{k(t)\langle D_x \rangle^{\frac{1}{\theta}}_{h} } P_\Lambda e^{-k(t) \langle D_x \rangle^{\frac{1}{\theta}}_{h}} $}

In order to apply the pseudodifferential calculus we need the condition $2\mu -1 < \theta$, which is fulfiled by taking $\mu$ close to $1$. We recall that $\theta < \frac{1}{q}$, so the $\xi-$order of $k(t)\langle \xi \rangle^{\frac{1}{\theta}}$ is higher than the $\xi-$order of $\Lambda$. We are going to consider $k(t)$ as in $\eqref{k(t)}$
for some positive $\rho_0 > 0$. We remark that
$$
e^{k(t)\langle D_x \rangle^{\frac{1}{\theta}}_{h} } D_t e^{-k(t) \langle D_x \rangle^{\frac{1}{\theta}}_{h}} = D_t - \frac1ik'(t)\langle D \rangle^{\frac{1}{\theta}}_{h},
\qquad \quad
e^{k(t)\langle D_x \rangle^{\frac{1}{\theta}}_{h} } a_3(t) D^{3}_{x} e^{-k(t) \langle D_x \rangle^{\frac{1}{\theta}}_{h}} = a_3(t)D^{3}_{x}.
$$
On the other hand, given $p(t,x,\xi) \in S^{m}_{\mu}(\R^2)$ we have
$$
e^{k(t)\langle D_x \rangle^{\frac{1}{\theta}}_{h} } p(t,x,D) e^{-k(t) \langle D_x \rangle^{\frac{1}{\theta}}_{h}} = q_1 (t,x,D)+ r_{\infty}(t,x,D)
$$
where $r_{\infty} (t,x,\xi)\in \mathcal{K}_{2\mu-1}$ and
$$
q_1(t,x,\xi) \sim \sum_{\gamma} \frac{1}{\gamma!} \partial^{\gamma}_{\xi} e^{k(t) \langle \xi \rangle^{\frac{1}{\theta}}_{h} } D^{\gamma}_{x}\{p(t,x,\xi) e^{-k(t) \langle \xi \rangle^{\frac{1}{\theta}}_{h}}\} \quad \text{in} \quad FS^{\infty}_{\mu;\frac{1}{\theta}}.
$$
Since 
$$
|\partial^{\alpha}_{\xi} \partial^{\beta}_{x} \{\partial^{\gamma}_{\xi} e^{k(t) \langle \xi \rangle^{\frac{1}{\theta}}_{h} } D^{\gamma}_{x}\{p(t,x,\xi) e^{-k(t) \langle \xi \rangle^{\frac{1}{\theta}}_{h}} \}\}| \leq C^{\alpha+\beta+2\gamma+1} (\alpha!\beta!\gamma!^{2})^{\mu} \langle \xi \rangle^{m-\gamma(1-\frac{1}{\theta}) -\alpha}_{h},
$$
we actually have $q_1 (t,x,\xi)\in S^{m}_{\mu}(\R^2)$.
Hence, we may write 
\begin{equation}\label{bullet}
e^{k(t)\langle D_x\rangle^{\frac{1}{\theta}}_{h} } p(t,x,D) e^{-k(t) \langle D_x \rangle^{\frac{1}{\theta}}_{h}} = p(t,x,D) + p_{k(t)}(t,x,D) + r_{\infty}(t,x,D),
\end{equation}
where $p_{k(t)}(t,x,\xi) \in S^{m-(1-\frac{1}{\theta})}_{\mu}$. Moreover, if in addition $p$ satisfies 
$$
|\partial^{\alpha}_{\xi}\partial^{\beta}_{x}p(t,x,\xi)| \leq t^{s}C^{\alpha+\beta+1} \langle \xi \rangle^{m-\alpha}_{h} \langle x \rangle^{-\sigma},
$$
then $p_{k(t)}$ satisfies
\begin{equation}\label{goomba}
	|\partial^{\alpha}_{\xi}\partial^{\beta}_{x}p_{k(t)}(t,x,\xi)| \leq t^{s}C^{\alpha+\beta+1} \langle \xi \rangle^{m-(1-\frac{1}{\theta})-\alpha}_{h} \langle x \rangle^{-\sigma}.
\end{equation}

Taking into account \eqref{bullet} and \eqref{goomba} we have
\begin{align*}
	iP_{k,\Lambda}(t,x,D) &= D_t -\frac1i k'(t) \langle D \rangle^{\frac{1}{\theta}}_{h} - D_t\Lambda (t,x,D)\\  
	&+ a_3(t) D^3_x + ia_3(t) op(3\xi^2 \partial_x \Lambda) 
	+ a_2(t,x)D^{2}_{x} + r_{2}(t,x,D)\\ 
	&+ a_1(t,x)D_{x} + r_{1}(t,x,D) + (r_{0} + r_{\infty})(t,x,D),
\end{align*}
where $r_{2}(t,x,\xi) \in SG^{2-(1-\frac{1}{\theta}), -\sigma_2}_{\mu}$, $r_{1}(t,x,\xi) \in SG^{1-(1-\frac{1}{\theta}), -\sigma_1}_{\mu}$, $(r_0 + r_{\infty})(t,x,\xi)$ is of order zero and
$$
|\partial^{\alpha}_{\xi} \partial^{\beta}_{x} r_{2}(t,x,\xi)| \leq t^{k} C^{\alpha+\beta+1} (\alpha!\beta!)^{\mu} \langle \xi \rangle^{2-(1-\frac{1}{\theta})-\alpha}_{h} \langle x \rangle^{-\sigma_2},$$
$$|\partial^{\alpha}_{\xi} \partial^{\beta}_{x} {r}_{1}(t,x,\xi)| \leq t^{k'} C^{\alpha+\beta+1} (\alpha!\beta!)^{\mu} \langle \xi \rangle^{1-(1-\frac{1}{\theta})-\alpha}_{h} \langle x \rangle^{-\sigma_1}.
$$
We point out that $r_{2}$ and $r_{1}$ have symbols depending on the whole transformation, that is, they depend on all the parameters $M_1, M_2, M_{e,1}, M_{e,2}, M_{\psi,1}, M_{\psi,2}$ and $\rho_0$. Moreover, notice that 
\beqs\label{h2}
|r_{2}(t,x,\xi)| &\leq& t^{k} C(\rho_0,\Lambda) h^{-(1-\frac{1}{\theta})} \langle \xi \rangle_h^{2}\langle x \rangle^{-\sigma_2},
\\
\label{h1}
|{r}_{1}(t,x,\xi)| &\leq& t^{k'} C(\rho_0,\Lambda)h^{-(1-\frac{1}{\theta})} \langle \xi \rangle_h \langle x \rangle^{-\sigma_1}.
\eeqs


\subsection{Lower bound estimates and energy estimate}
In this section we shall use sharp-G{\aa}rding and Fefferman-Phong inequalities in order to get a priori energy estimates for the Cauchy problem \eqref{CPlambda}. We first write 
\begin{align*}
	iP_{k,\Lambda}(t,x,D) &= \partial_t + \rho_0 \langle D_x \rangle^{\frac{1}{\theta}}_{h} - \partial_t \Lambda_{e} (t,D)- \partial_{t} \Lambda_{\psi} (t,D)\\
	&+ ia_3(t) D^3_x - a_3(t) op( 3\xi^2 \partial_x \Lambda )\\
	&+ ia_2(t,x)D^{2}_{x} + ir_2(t,x,D) \\
	&+ ia_1(t,x)D_{x} + ir_1 (t,x,D)+ r_{0}(t,x,D).
\end{align*}
Since we want to apply Fefferman-Phong inequality at level 2, we split the imaginary part of the terms of order $2$ into its hermitian and anti-hermitian parts:
$$Im\left(ia_2(t,x)D^{2}_{x} + ir_2(t,x,D)\right) = Rea_2(t,x) D_x^2+Re r_2(t,x,D)=A(t,x,D),$$
and we rewrite the operator $iA$ as
\beqsn
iA&=&\frac{iA+(iA)^*}2+\frac{iA-(iA)^*}2
=\frac{1}{2} op\left(\sum_{\gamma \geq 1} \frac{1}{\gamma!} \partial^{\gamma}_{\xi}D^{\gamma}_{x} \overline{iA}\right) +A_2 + r
\eeqsn
with $Re \langle A_2(t,x,D) v, v \rangle = 0$ for all $v\in\mathscr{S}(\R)$, $r(t,x,\xi)\in\mathcal{K}_{2\mu-1},$ so we get
\begin{align*}
	ia_2(t,x)D^{2}_{x} + ir_2(t,x,D) &=- Im\, a_2(t,x) D_x^2- Im\, r_2(t,x,D)+iA(t,x,D)
\\&= - Im\, a_2(t,x) D_x^2 - Im\, r_2(t,x,D) + A_{2} (t,x,D)
\\
&- (\partial_x Re a_2(t,x)) D_x - \frac{i}{2} op\left(\sum_{\gamma \geq 1} \frac{1}{\gamma!} \partial^{\gamma}_{\xi}D^{\gamma}_{x} Re\, r_2\right)  + r_0(t,x,D)\\
	& = - Im\, a_2(t,x) D_x^2 - Im\, r_2 + A_{2} 
	+ \tilde{r}_2 + r_0,
\end{align*}
where $\tilde{r}_{2}(t,x,\xi) = - (\partial_x\Re a_2) \xi - \frac{i}{2} \sum_{\gamma \geq 1} \frac{1}{\gamma!} \partial^{\gamma}_{\xi}D^{\gamma}_{x} Re\, r_2(t,x,\xi)$,
\beqs\label{h3}
|\tilde{r}_2| \lesssim t^{k} \langle \xi \rangle^{1}_h \langle x \rangle^{-\sigma_2}\lesssim t^{k} \langle \xi \rangle_h^2h^{-1} \langle x \rangle^{-\sigma_2}
\eeqs
and $r_0$ is of order zero. We thus come to
\begin{align*}
	iP_{k,\Lambda}(t,x,D)&=e^{k(t) \langle D \rangle^{\frac{1}{\theta}}_{h}} e^{\Lambda} iP(t,x,D) \{e^{\Lambda}\}^{-1}  e^{-k(t) \langle D_x \rangle^{\frac{1}{\theta}}_{h} }
\\
&= \partial_t + \rho_0 \langle D_x \rangle^{\frac{1}{\theta}}_{h} - \partial_t \Lambda_{e} (t,D)- \partial_{t} \Lambda_{\psi}(t,D)\\
	&+ ia_3(t) D^3_x  + A_{2}(t,x,D) \\
	&- Im\, a_2(t,x) D_x^2 - Im\, r_2 (t,x,D) + \tilde{r}_2(t,x,D) - a_3(t) op(3\xi^2 \partial_x \Lambda) \\
	&+ ia_1(t,x)D_{x} + ir_1(t,x,D) 
	+ r_{0}(t,x,D).
\end{align*}
In the following computations we shall assume $|\xi| > 2h$ (in the complementary region all symbols are of order zero with seminorms depending on the parameter $h$). We notice that by definition \eqref{defl}
$$\partial_x \Lambda=\partial_x \lambda_{e,2}+\partial_x \lambda_{e,1}$$ and we write

\begin{align*}
	- a_3(t) 3\xi^2 \partial_x \lambda_{e,2} (t,x,\xi)&= |a_3(t)|t^{-(\ell-k)}3M_2 \xi^2 \langle x \rangle^{-\sigma_2} (1-\chi)\left( t\langle \xi \rangle^{\frac{2-q_2}{k+1}}_{h} \right) \psi \left( \frac{\langle x\rangle}{\langle \xi \rangle^{\frac{2-q_2}{\sigma_2}}_{h}} \right) \\
	&= |a_3(t)|t^{-(\ell-k)}3M_2 \xi^2 \langle x \rangle^{-\sigma_2} (1-\chi)\left( t\langle \xi \rangle^{\frac{2-q_2}{k+1}}_{h} \right) \\
	&- |a_3(t)|t^{-(\ell-k)}3M_2 \xi^2 \langle x \rangle^{-\sigma_2} (1-\chi)\left( t\langle \xi \rangle^{\frac{2-q_2}{k+1}}_{h} \right) (1-\psi) \left( \frac{\langle x\rangle}{\langle \xi \rangle^{\frac{2-q_2}{\sigma_2}}_{h}} \right) \\
	&= |a_3(t)|t^{-(\ell-k)}3M_2 \xi^2 \langle x \rangle^{-\sigma_2} (1-\chi)\left( t\langle \xi \rangle^{\frac{2-q_2}{k+1}}_{h} \right) + r_{2,\psi}(t,x,\xi),
\end{align*}
where 
\beqs\label{j2}|r_{2,\psi}(t,x,\xi)|\lesssim t^k\langle\xi\rangle_h^2\langle x\rangle^{-\sigma_2}(1-\psi) \left( \frac{\langle x\rangle}{\langle \xi \rangle^{\frac{2-q_2}{\sigma_2}}_{h}} \right)\lesssim t^k\langle\xi\rangle_h^2\langle \xi\rangle_h^{-2+q_2}=t^k\langle\xi\rangle_h^{q_2}\eeqs is of order $q_2$ without any decay. Similarly we get
\begin{align*}
	- a_3(t) 3\xi^2 \partial_x \lambda_{e,1}(t,x,\xi) &=|a_3(t)| t^{-(\ell-k')}3M_1 \xi^2 \langle \xi \rangle^{-1}_{h} \langle x \rangle^{-\sigma_1} (1-\chi)\left( t\langle \xi \rangle^{\frac{1-q_1}{k'+1}}_{h} \right) \psi \left( \frac{\langle x\rangle}{\langle \xi \rangle^{\frac{1-q_1}{\sigma_1}}_{h}} \right) \\
	&= |a_3(t)| t^{-(\ell-k')}3M_1 \xi^2 \langle \xi \rangle^{-1}_{h} \langle x \rangle^{-\sigma_1} (1-\chi)\left( t\langle \xi \rangle^{\frac{1-q_1}{k'+1}}_{h} \right) \\
	&- |a_3(t)| t^{-(\ell-k')}3M_1 \xi^2 \langle \xi \rangle^{-1}_{h} \langle x \rangle^{-\sigma_1} (1-\chi)\left( t\langle \xi \rangle^{\frac{1-q_1}{k'+1}}_{h} \right) (1-\psi) \left( \frac{\langle x\rangle}{\langle \xi \rangle^{\frac{1-q_1}{\sigma_1}}_{h}} \right) \\
	&= |a_3(t)| t^{-(\ell-k')}3M_1 \xi^2 \langle \xi \rangle^{-1}_{h} \langle x \rangle^{-\sigma_1} (1-\chi)\left( t\langle \xi \rangle^{\frac{1-q_1}{k'+1}}_{h} \right) + r_{1,\psi}(t,x,\xi),
\end{align*}
where 
\beqs\label{j1}r_{1,\psi}(t,x,\xi)\leq t^{k'}\langle\xi\rangle_h^{q_1}\eeqs is of order $q_1$ without any decay. Since $q=\max\{q_2,q_1\}$, we can say that the operator $r_{2,\psi}(t,x,D)+r_{1,\psi}(t,x,D)$ has order $q$, and we shall use the term $\rho_0 \langle D_x\rangle^{\frac{1}{\theta}}_{h}$ to control it. Then we may write

\beqsn
	iP_{k,\Lambda}(t,x,D)
&=& \partial_t  
	  + ia_3(t)D^{3}_{x}+ A_2- \partial_t \Lambda_{e}(t,D) \\
&+&op\left\{\left(-Im\,a_2 \xi^2 - Im\, r_2+\tilde r_2 + |a_3(t)| t^{-(\ell-k)}3M_2 \xi^2 \langle x \rangle^{-\sigma_2} \right)(1-\chi)\left( t\langle \xi \rangle^{\frac{2-q_2}{k+1}}_{h} \right)\right\} \\
&+&op\left \{(- Im\, a_2 \xi^2 - Im\, r_2+\tilde r_2) \chi \left( t \langle \xi \rangle^{\frac{2-q_2}{k+1}}_{h}\right) \right\}- \partial_t \Lambda_{\psi,2} (t,D)
\\
&+&
op\left\{\left(i a_1 \xi   + ir_1 + |a_3(t)| t^{-(\ell-k')}3M_1 \xi^2 \langle \xi \rangle^{-1}_{h} \langle x \rangle^{-\sigma_1}\right) (1-\chi)\left( t\langle \xi \rangle^{\frac{1-q_1}{k'+1}}_{h} \right)\right\}\\
	&+& op\left\{(i a_1 \xi + ir_1) \chi \left( t \langle \xi \rangle^{\frac{1-q_1}{k'+1}}_{h}\right) \right\}- \partial_t \Lambda_{\psi,1}(t,D) \\
&
	+&\rho_0 \langle D_x \rangle^{\frac{1}{\theta}}_{h} + r_{2,\psi}(t,x,D) + r_{1,\psi} (t,x,D)+ r_0(t,x,D).
\eeqsn

In the following we shall set up the parameters in order to make positive the symbols of the lower order terms of $iP_{k, \Lambda}$. We begin with the terms of order $2$. In the evolution zone we have, using \eqref{a3}, \eqref{a2}, \eqref{h2}, \eqref{h3}
\begin{align*}
&\left(-Im\,a_2 \xi^2 - Im\, r_2+\tilde r_2 + |a_3(t)| t^{-(\ell-k)}3M_2 \xi^2 \langle x \rangle^{-\sigma_2} \right)(1-\chi)\left( t\langle \xi \rangle^{\frac{2-q_2}{k+1}}_{h} \right)
\\
	&\geq t^{k} \xi^2 \langle x \rangle^{-\sigma_2}\{- C_{a_2} - C(\Lambda, \rho_0)h^{-(1-\frac{1}{\theta})}-C(\Lambda, \rho_0)h^{-1}+3c_{a_3}M_2\} (1-\chi) \left( t \langle \xi \rangle^{\frac{2-q_2}{k+1}}_{h}\right),
\end{align*}
whereas in the pseudodifferential zone, using \eqref{a2}, \eqref{h2},\eqref{LP2},  \eqref{h3} we get 
\begin{align*}
	&-\partial_{t}  \Lambda_{\psi, 2} + \{-Im\,a_2 \xi^2 - Im\, r_2+\tilde r_2\} \chi \left( t \langle \xi \rangle^{\frac{2-q_2}{k+1}}_{h}\right) \\
	&= \{M_{\psi,2} t^{k }\langle \xi \rangle^{2}_{h } -Im\,a_2 \xi^2 - Im\, r_2 +\tilde r_2\} \chi \left( t \langle \xi \rangle^{\frac{2-q_2}{k+1}}_{h}\right)\\
	&\geq t^{k} \langle \xi \rangle^{2}_{h} \{ M_{\psi,2} - C_{a_2} - C(\Lambda, \rho_0) h^{-(1-\frac{1}{\theta})} -C(\Lambda, \rho_0)h^{-1}\}\chi \left( t \langle \xi \rangle^{\frac{2-q_2}{k+1}}_{h}\right).
\end{align*}
Choosing $M_2$ and $M_{\psi,2}$ large (depending on $a_2$), precisely 
$$M_2>C_{a_2}/(3C_{a_3})\ {\rm and} \ M_{\psi,2}> C_{a_2},$$ we make the symbol of level two non-negative as long as we take $h$ (depending on all the parameters) large enough. Hence, we can apply the Fefferman-Phong inequality to the terms of level 2, say $iP_{k,\Lambda,2}$ for brevity's sake, and we get 
$$
Re\langle iP_{k,\Lambda,2}(t,x,D)v,v\rangle_{L^2}\geq -C\|v\|_{L^2}, \quad\forall v\in\mathscr S(\R).
$$

For the evolution zone of terms of order $1$ we observe that, by \eqref{a3}, \eqref{a2}, \eqref{a1}, \eqref{h1},
\begin{align*}
&\left(i a_1 \xi + ir_1 + |a_3(t)| t^{-(\ell-k')}3M_1 \xi^2 \langle \xi \rangle^{-1}_{h} \langle x \rangle^{-\sigma_1}\right) (1-\chi)\left( t\langle \xi \rangle^{\frac{1-q_1}{k'+1}}_{h} \right)
 \\
	&\geq t^{k'}|\xi|\langle x \rangle^{-\sigma_1} \left\{ - C_{a_1} - C(\Lambda, \rho_0) h^{-(1-\frac{1}{\theta})} + \frac{3c_{a_3}M_1}{2} \right\}
	(1-\chi) \left( t \langle \xi \rangle^{\frac{1-q_1}{k'+1}}_{h}\right),
\end{align*}
whereas for the pseudodifferential zone, using \eqref{a2}, \eqref{a1}, \eqref{LP1}, \eqref{h1},  we have 
\begin{align*}
	-\partial_{t}  \Lambda_{\psi, 1} &+ \{i a_1 \xi + ir_1\} \chi \left( t \langle \xi \rangle^{\frac{1-q_1}{k'+1}}_{h}\right) \\
	&= \{M_{\psi,1} t^{k'}\langle \xi \rangle_{h } + i a_1 \xi - \partial_x Re\, a_2 \xi  + ir_1 \} \chi \left( t \langle \xi \rangle^{\frac{1-q_1}{k'+1}}_{h}\right)\\
	&\geq t^{k'} \langle \xi \rangle_{h} \{ M_{\psi,1} - C_{a_1} - C(\Lambda, \rho_0) h^{-(1-\frac{1}{\theta})} \}\chi \left( t \langle \xi \rangle^{\frac{1-q_1}{k'+1}}_{h}\right).
\end{align*}
Againg we can turn the symbol of level 1 into a non-negative symbol by choosing $M_{1}$, $M_{\psi,1}$ large, precisely 
$$
M_1>2C_{a_1}/(3c_{a_3}) \ {\rm and}\ M_{\psi,1} >C_{a_1},
$$
and $h$ large enough depending on all the parameters. Hence, we can apply the sharp-Garding inequality to the terms of level 1, say $iP_{k,\Lambda,1}$ for brevity's sake, and we get 
$$
Re\langle iP_{k,\Lambda,1}(t,x,D)v,v\rangle_{L^2}\geq -C\|v\|_{L^2},\quad  \forall v\in\mathscr S(\R).
$$

Now, to deal with $-\partial_t \Lambda_{e}(t,\xi)$ we note that
\begin{small}
	\begin{align*}
		-\partial_t\Lambda_e &= t^{-(\ell-k) - 1} \left\{ (\ell-k) \lambda_2 (1-\chi)\left( t\langle\xi\rangle_h^{\frac{2-q_2}{k+1}} \right) 
		+\lambda_2 \chi'\left( t\langle\xi\rangle_h^{\frac{2-q_2}{k+1}} \right)t\langle\xi\rangle_h^{\frac{2-q_2}{k+1}}
		+ M_{e,2} \langle \xi \rangle^{(2-q_2)\frac{1-\sigma_2}{\sigma_2}}_{h}  \right\} \\
		&\times (1-\chi)\left( 2t\langle\xi\rangle_h^{\frac{2-q_2}{k+1}} \right) \\
		&+ t^{-(\ell-k') - 1} \left\{ (\ell-k') \lambda_1 (1-\chi)\left( t\langle\xi\rangle_h^{\frac{1-q_1}{k'+1}} \right) 
		+\lambda_1 \chi'\left( t\langle\xi\rangle_h^{\frac{1-q_1}{k'+1}} \right)t\langle\xi\rangle_h^{\frac{1-q_1}{k'+1}}
		+ M_{e,1} \langle \xi \rangle^{(1-q_1)\frac{1-\sigma_1}{\sigma_1}}_{h}  \right\} \\
		&\times (1-\chi)\left( 2t\langle\xi\rangle_h^{\frac{1-q_1}{k'+1}} \right),
	\end{align*}
\end{small}
so, looking at \eqref{estl2} and \eqref{estl1}, we may choose $M_{e,2}$ (depending on $M_2$) and $M_{e,1}$ (depending on $M_1$) to make $-\partial_t \Lambda_{e}$ a non-negative symbol as well. Since the order of $\partial_t \Lambda$ is at most two, see \eqref{star2} and \eqref{star1}, we can apply Fefferman-Phong inequality to $-\partial_t\Lambda_e(t,D)$ and we get
$$Re\langle-\partial_t\Lambda_e(t,D)v,v\rangle_{L^2}\geq -C\|v\|_{L^2},\quad  \forall v\in \mathscr{S}(\R).$$

It remains to consider $\rho_0 \langle D_x\rangle^{\frac{1}{\theta}}_{h} + r_{2,\psi} (t,x,D)+ r_{1,\psi}(t,x,D)$. We have by \eqref{j2} and \eqref{j1}
\begin{align*}
\rho_0 \langle \xi \rangle^{\frac{1}{\theta}}_{h} + r_{2,\psi}(t,x,\xi) + r_{1,\psi} (t,x,\xi)&\geq \rho_0 \langle \xi \rangle^{\frac{1}{\theta}}_{h} - c(M_2) \langle \xi \rangle^{q_2}_{h} - c(M_1) \langle \xi \rangle^{q_1}_{h} \\
&\geq \langle \xi \rangle^{\frac{1}{\theta}}_{h} \{ \rho_0 - c(M_2)h^{-(\frac{1}{\theta}-q_2)} - c(M_1)h^{-(\frac{1}{\theta}-q_1)}\}.
\end{align*}
Hence, as long as we choose $h$ large enough, we have $\rho_0 \langle \xi \rangle^{\frac{1}{\theta}}_{h} + r_{2,\psi} (t,x,\xi)+ r_{1,\psi}(t,x,\xi)\geq 0$ for any fixed $\rho_0 > 0$, and we can apply the sharp-G{\aa}rding inequality to this term.

Summing up, applying Fefferman-Phong inequality to the real terms of order $2$ and sharp-G{\aa}rding inequality to the terms of order at most $1$ we can write 
$$
iP_{k,\Lambda} (t,x,D)= \partial_t +ia_3(t)D^{3}_{x} + A_2(t,x,D) + p(t,x,D) + r_0(t,x,D),
$$
where $r_0$ is of order zero, $Re\, \langle A_2(t,x,D) v,v \rangle = 0$ and  $Re\, \langle p(t,x,D)v, v \rangle \geq -C \|v\|_{L^2}$ for every $v\in\mathscr S(\R)$.
Therefore 
\begin{align*}
	\partial_t \|v(t)\|^{2}_{L^2} &= 2Re\, \langle \partial_t v, v \rangle \\
	& = 2Re\, \langle \{iP_{k,\Lambda}(t,x,D) - ia_3(t)D^{3}_{x} - A_2(t,x,D) - p (t,x,D)- r_0(t,x,D)\} v, v \rangle \\
	&\leq \|P_{k,\Lambda}v\|^{2}_{L^2} + C_h \|v\|^{2}_{L^2}.
\end{align*}
Gronwall inequality then gives 
\begin{equation*}
	\|v(t)\|^{2}_{L^2} \leq C_h \left\{ \|v(0)\|^{2}_{L^2} + \int_0^t \|P_{k,\Lambda}v(\tau)\|^2_{L^2} d\tau \right\}.
\end{equation*}

To obtain the above inequality for general $H^{m}(\R)$ norms ($m\geq 0$) in the place of $L^2$ norms, we remark that the operator 
$$
P_{m}(t,x,D) = \langle D_x \rangle^{m}_{h} P(t,x,D) \langle D_x \rangle^{-m}_{h} = P (t,x,D)+ b_{m,2}(t,x,D) + b_{m,0}(t,x,D),
$$ 
where $b_{m,0}(t,x,\xi) \in S^{0}_{\mu}(\R^2)$ and 
$$
|\partial^{\alpha}_{\xi}\partial^{\beta}_{x}b_{m,2}(t,x,\xi)| \leq t^{k}C_m C^{\alpha+\beta} \alpha!\beta! \langle \xi \rangle^{1-\alpha}_{h} \langle x \rangle^{-\sigma_2}.
$$
Then we can treat $b_{m,2}$ as a term of order $2$, whose order is strictly less than $2$ and the decay is exactly the decay of the coefficient $a_2(t,x)$. In this way, repeating the conjugation argument for the operator $P_m$ and taking $h > h_0(\Lambda, \rho_0, m)$, we have 
\begin{equation*}
	\|v(t)\|^{2}_{H^m} \leq C_{h,m} \left\{ \|v(0)\|^{2}_{H^m} + \int_0^t \|P_{k,\Lambda}v(\tau)\|^2_{H^m} d\tau \right\}.
\end{equation*}

Let us summarize what we have proven in the following proposition.

\begin{proposition}
	Let $\rho_0 > 0$ and $m \geq 0$. There exist positive constants $M_2, M_1, M_{e,2}, M_{e,1}$, $M_{\psi,2}, M_{\psi, 1}$ and a parameter $h_0 = h_0(M_2, M_1, M_{e,2}, M_{e,1}, M_{\psi,2}, M_{\psi, 1}, \rho_0, m)$ such that for all $h \geq h_0$ the following a priori energy estimate holds
	\begin{equation}
		\|v(t)\|^{2}_{H^m} \leq C_{h,m} \left\{ \|v(0)\|^{2}_{H^m} + \int_0^t \|P_{k,\Lambda}v(\tau)\|^2_{H^m} d\tau \right\}
	\end{equation}
	for every $v \in C^{1}([0,T]; H^{m}(\R)) \cap C([0,T];H^{m+3}(\R))$.
\end{proposition}

By standard energy method arguments the above propostion implies that the Cauchy problem \eqref{CPlambda} for the operator $P_{k,\Lambda}$ is well-posed in $L^2(\R)$. As a consequence, we get well-posedness in $H^{\infty}_{\theta}(\R)$ for \eqref{genCP}, provided that $\theta \in (1,\frac{1}{q})$. This ends the proof of Theorem \ref{main} in the case $\sigma_2, \sigma_1 \in (0,1)$.

\subsection{Theorem \ref{main}, $\sigma_2 \in (0,1)$ and $\sigma_1 \geq 1$}

 The proof is still based on a change of variable $v = e^{k(t) \langle D_x \rangle^{\frac{1}{\theta}}_{h} } e^{\Lambda}(t,x,D)u$ of the form \eqref{ch}. We will only point out the differences in the definition of $\Lambda$, since all the rest of the proof can be obtained in the same way.

Since $\sigma_2 \in (0,1)$ and $\sigma_1 \geq 1$, we keep the symbols $\Lambda_2, \Lambda_{\psi, 1}, k(t)\langle\xi \rangle^{\frac{1}{\theta}}_{h}$ exactly as in \eqref{L2}, \eqref{lE2}, \eqref{LE2}, \eqref{LP2}, \eqref{LP1}, \eqref{k(t)} (hence the index $q_2$ remains unchanged). On the other hand, we note now that $\lambda_1$, still defined as in \eqref{L1}, now satisfies for $\alpha,\beta \in \N_0$

\beqs\label{panda_1}
|\partial^{\alpha}_{\xi} \partial^{\beta}_{x} \lambda_1(x,\xi)| &\leq&  M_1
	\begin{cases}
	C^{\alpha+\beta+1} \alpha!^{\mu} \beta!^{\mu} \langle \xi \rangle^{-1-\alpha}_{h} \log \langle x \rangle \langle x \rangle^{-\beta} \\
	C^{\alpha+\beta+1} \alpha!^{\mu} \beta!^{\mu} \log \langle \xi \rangle_h \langle \xi \rangle^{ - 1 -\alpha}_{h} \langle x \rangle^{-\beta}, 
	\end{cases}
	\quad \sigma_1 = 1,
	\\\label{panda11}
|\partial^{\alpha}_{\xi} \partial^{\beta}_{x} \lambda_1(x,\xi)| &\leq& M_1
	C^{\alpha+\beta+1} \alpha!^{\mu} \beta!^{\mu} \langle \xi \rangle^{-1-\alpha}_{h} \langle x \rangle^{-\beta}, \qquad\qquad\sigma_1 > 1.
\eeqs

Therefore $\lambda_{e,1}$ (still defined as in \eqref{lE1}) now satisfies 
\beqsn
|\partial^{\alpha}_{\xi} \partial^{\beta}_{x} \lambda_{e,1}(t,x,\xi)| \leq M_1
\begin{cases}
	C^{\alpha+\beta+1} \alpha!^{\mu} \beta!^{\mu} \langle \xi \rangle^{(1-q_1) \frac{\ell-k'}{k'+1}-1-\alpha}_{h} \log \langle x \rangle \langle x\rangle^{-\beta}, \,\,\,\,\, \sigma_1 = 1, \\
	C^{\alpha+\beta+1} \alpha!^{\mu} \beta!^{\mu} \log \langle \xi \rangle_{h}\langle \xi \rangle^{(1-q_1) \frac{\ell-k'}{k'+1}-1-\alpha}_{h} \langle x \rangle^{-\beta}, \,\,\, \sigma_1 = 1, \\
	C^{\alpha+\beta+1} \alpha!^{\mu} \beta!^{\mu} \langle \xi \rangle^{(1-q_1) \frac{\ell-k'}{k'+1}-1-\alpha}_{h} \langle x\rangle^{-\beta}, \,\,\,\,\quad \qquad \sigma_1 > 1
\end{cases}
\eeqsn
so we define 
$$
\Lambda_{e,1}(t,\xi) = -M_{e,1} 
\begin{cases}
	 \log \langle \xi \rangle_h \langle \xi \rangle^{-1}  \ds\int_0^t\tau^{k'-\ell-1} (1-\chi)\left( 2\tau\langle\xi\rangle_h^{\frac{1-q_1}{k'+1}} \right) d\tau, \quad \sigma_1 = 1 \\
	 \langle \xi \rangle^{-1}  \ds\int_0^t\tau^{k'-\ell-1} (1-\chi)\left( 2\tau\langle\xi\rangle_h^{\frac{1-q_1}{k'+1}} \right) d\tau, \quad \sigma_1 > 1.
\end{cases}
$$
Since $\Lambda_{\psi,1}$ and $\lambda_{e,1}+\Lambda_{e,1}$ should have the same order, the new condition on $q_1$ is
$$
q_1 = (1-q_1)\frac{\ell-k'}{k' + 1} - 1 \iff q_1 = \frac{\ell-2k'-1}{\ell+1},
$$
see the second line in \eqref{qq1}.
We remark that $q_1$ has an influence on the Gevrey index $\theta$ only if $\ell > 2k' + 1$.


\subsection{Theorem \ref{main}, $\sigma_1 \in (0,1)$ and $\sigma_2 \geq 1$}

In this case we keep the symbols $\Lambda_1, \Lambda_{\psi, 2}, k(t)\langle\xi \rangle^{\frac{1}{\theta}}_{h}$ exactly as in 
\eqref{L1}, \eqref{lE1}, \eqref{LE1}, \eqref{LP1}, \eqref{LP2}, \eqref{k(t)} (hence the index $q_1$ remains unchanged). On the other hand, we note now that $\lambda_2$, still defined as in \eqref{L2}, now satisfies for $\alpha,\beta \in \N_0$

\beqs\label{panda_2}
|\partial^{\alpha}_{\xi} \partial^{\beta}_{x} \lambda_2(x,\xi)| &\leq& M_2
\begin{cases}
	C^{\alpha+\beta+1} \alpha!^{\mu} \beta!^{\mu} \langle \xi \rangle^{-\alpha}_{h} \log \langle x \rangle \langle x \rangle^{-\beta}, \\
	C^{\alpha+\beta+1} \alpha!^{\mu} \beta!^{\mu} \log \langle \xi \rangle_h \langle \xi \rangle^{ -\alpha}_{h} \langle x \rangle^{-\beta}, 
\end{cases}
\quad \sigma_2 = 1,
\eeqs
\beqs\label{2panda2}
|\partial^{\alpha}_{\xi} \partial^{\beta}_{x} \lambda_2(x,\xi)| &\leq& M_2
	C^{\alpha+\beta+1} \alpha!^{\mu} \beta!^{\mu} \langle \xi \rangle^{-\alpha}_{h} \langle x \rangle^{-\beta}, \qquad\qquad\sigma_2 > 1.
\eeqs
Therefore $\lambda_{e,2}$ (still defined as in \eqref{lE2}) now satisfies
$$
|\partial^{\alpha}_{\xi} \partial^{\beta}_{x} \lambda_{e,2}(t,x,\xi)| \leq M_2
\begin{cases} 
	C^{\alpha+\beta+1} \alpha!^{\mu} \beta!^{\mu} \langle \xi \rangle^{(2-q_2) \frac{\ell-k}{k+1}-\alpha}_{h} \log \langle x \rangle \langle x\rangle^{-\beta},\quad \sigma_2=1 \\
	C^{\alpha+\beta+1} \alpha!^{\mu} \beta!^{\mu} \log \langle \xi \rangle_{h}\langle \xi \rangle^{(2-q_2) \frac{\ell-k}{k+1}-\alpha}_{h} \langle x \rangle^{-\beta}, \quad \sigma_2>1
\end{cases}
$$
and so we define 
$$
\Lambda_{e,2}(t,\xi) = -M_{e,2}
\begin{cases}
	 \log \langle \xi \rangle_h  \ds \int_0^t \tau^{k-\ell-1} (1-\chi)\left(2\tau\langle\xi\rangle_h^{\frac{2-q_2}{k+1}} \right) d\tau, \quad \sigma_2 = 1 \\
	 \ds\int_0^t \tau^{k-\ell-1} (1-\chi)\left(2\tau\langle\xi\rangle_h^{\frac{2-q_2}{k+1}} \right) d\tau, \quad \sigma_2 > 1.
\end{cases}
$$
Since $\Lambda_{\psi,2}$ and $\lambda_{e,2}+\Lambda_{e,2}$ should have the same order, the new condition on $q_2$ is
$$
q_2 = (2-q_2)\frac{\ell-k}{k + 1} \iff q_2 = 2\frac{\ell-k}{\ell+1},
$$
see the second line in \eqref{qq2}.
Note that $0 < q_2 < 1$, since $\ell > k$ and $\ell < 2k  +1$ by hypothesis.


\subsection{Theorem \ref{main}, $\sigma_2, \sigma_1 \geq 1$.}

In this last (and most simple) case the symbols $\lambda_2$ and $\lambda_1$, given in $\eqref{L2}$ and $\eqref{L1}$, satisfy \eqref{panda_2}/\eqref{2panda2} and \eqref{panda_1}/\eqref{panda11}, respectively. Therefore the indices $q_2$ and $q_1$ in this case are 
$$
q_2 = 2\frac{\ell-k}{\ell+1}, \quad q_1 = \frac{\ell-2k'-1}{\ell+1}.
$$
Since $\ell > k$ and $\ell < 2k  +1$ we have $0 < q_2 < 1$, while $q_1 > 0$ only if $\ell > 2k'+1$.

\begin{remark}
		Notice that in the limit case $\sigma_2=1$ and $\ell=k$ (no slow decay and no degeneration gap in time at level 2) we have $q_2=0$, and the definition of $\lambda_2$ in \eqref{L2} turns out to be the one used in literature to treat the non-degenerate $H^\infty$ case, see \cite{KB}, \cite{CRJEECT}; similarly, if $\sigma_1=1$ and $\ell=k'$(no slow decay and no degeneration gap in time at level 1) the definition of $\lambda_1$ in \eqref{L1} turns out to be the one used in literature to treat the non-degenerate $H^\infty$ case, see \cite{ascanelli_chiara_zhanghirati_well_posedness_of_cauchy_problem_for_p_evo_equations}.
	\end{remark}


\section {Proof of Theorem \ref{main2}}

The proof of Theorem \ref{main2} is still based on the change of variable $v = e^{k(t) \langle D_x \rangle^{\frac{1}{\theta}}_{h} } e^{\Lambda}(t,x,D)u$. The advantage here is that we have no degeneracy gap in the term of order $2$ ($\ell \leq k$), so there's no need to introduce a pseudodifferential zone and an evolution zone at level 2. The function $k(t)$ is defined as in \eqref{k(t)} whereas the  symbol $\Lambda$ is defined as 
$$
\Lambda = \lambda_2 + \lambda_{e,1} + \Lambda_{e,1} + \Lambda_{\psi,1}, 
$$
where $\lambda_2$, $\lambda_{e,1}$, $\Lambda_{e,1}$, $\Lambda_{\psi, 1}$ are defined as in $\eqref{L2}, \eqref{lE1}, \eqref{LE1}, \eqref{LP1}$. Therefore the index $q_1$ does not change, it's the same as in the proof of Theorem \ref{main}.

\bigskip

\noindent \textbf{Case $\sigma_2 \in (\frac{1}{2},1)$}. Here $\lambda_2$ satisfies \eqref{estl2}. Therefore, the index $q_2$ here should satisfy 
$$
q_2 = (2-q_2)\frac{1-\sigma_2}{\sigma_2} \iff q_2 = 2(1-\sigma_2).
$$
Note that $q_2 < 1 \iff \sigma_2 > \frac{1}{2}$. Well-posedness in $H^\infty_\theta$ follows for $\theta<\min\{\frac1{2(1-\sigma_2)}\frac1{q_1}\}.$ 

\bigskip

\noindent \textbf{Case $\sigma_2 = 1$}. Here we have 
$$|\partial^{\alpha}_{\xi} \partial^{\beta}_{x} \lambda_2(x,\xi)| \leq M_2 C^{\alpha+1} \alpha!^{\mu}\beta!^{\mu} \log \langle \xi \rangle_h \langle \xi \rangle^{- \alpha}_{h} \langle x \rangle^{-\beta}, \quad \quad  \beta, \alpha \in \N_0.$$
In this way $e^{\lambda_2}$ is a symbol of finite order. Therefore:
\begin{itemize} 
\item[-] if $q_1 > 0$, then $e^\Lambda$ is a transformation of infinite order; using the change of variable $v = e^{k(t) \langle D_x \rangle^{\frac{1}{\theta}}_{h} } e^{\Lambda}u$ we get well-posedness in $H^{\infty}_{\theta}$ for all $1 < \theta <\frac{1}{q_1}$
\item[-] if $q_1 \leq 0$ we get well-posedness in $H^{\infty}(\R)$; indeed, in this last case the operator $e^\Lambda$ is of finite order, the transformation $e^{k(t) \langle D_x \rangle^{\frac{1}{\theta}}_{h}}$ is no more needed, and the change of variable $v = e^{\Lambda}u$ is enough to move the Cauchy problem \eqref{genCP} into an equivalent Cauchy problem well-posed in $L^2(\R)$ and $H^m(\R)$ for every $m$. Now, if $f(t),g$ are in $H^\infty(\R)$, then $e^\Lambda f(t), e^\Lambda g \in H^m(\R)$ for every $m$, so there exists a unique $v$ solution in $H^m(\R)$ of the auxiliary Cauchy problem. The solution of the original Cauchy problem is $u=(e^\Lambda)^{-1}v\in H^{m-\delta}(\R)$ for all $m$, so it is in $H^\infty(\R)$.
\end{itemize}

\bigskip

\noindent \textbf{Case $\sigma_2 > 1$}. Here we have 
$$|\partial^{\alpha}_{\xi} \partial^{\beta}_{x} \lambda_2(x,\xi)| \leq M_2 
	C^{\alpha+1} \alpha!^{\mu}\beta!^{\mu} \langle \xi \rangle^{- \alpha}_{h} \langle x \rangle^{-\beta}, \quad \quad  \beta, \alpha \in \N_0, $$
and in this way $e^{\lambda_2}$ is a symbol of order zero. Therefore if $q_1 > 0$ we get well-posedness in $H^{\infty}_{\theta}$ for all $1 < \theta <\frac{1}{q_1}$ by the change of variable $v = e^{k(t) \langle D_x \rangle^{\frac{1}{\theta}}_{h} } e^{\Lambda}u$, while if $q_1 \leq 0$ we get well-posedness in $L^{2}(\R)$ by the change of variable $v = e^{\Lambda}u$ of order zero.


\section{Proof of Theorem \ref{main3}.}

In this case we do not have to worry with the degeneracy gap in time, that is we have $t^{k} \leq t^{\ell}$ and $t^{k'} \leq t^{\ell}$ as $t \to 0^{+}$. No separation in pseudodifferential and evolution zones is needed, and the change of variable here is simply defined by $v = e^{k(t) \langle D_x\rangle^{\frac{1}{\theta}}}e^{\Lambda}(t,x,D)u$ where 
$$
k(t) =\begin{cases} 0, \quad \sigma_2 \geq 1 \,\, \text{and}\,\, \sigma_1 \geq \frac{1}{2}
\\
\eqref{k(t)}, \quad {\rm otherwise,}
\end{cases}
$$
whereas
$$
\Lambda = \lambda_2 + \lambda_1,
$$
with $\lambda_2$ given by $\eqref{L2}$ and $\lambda_1$ given by $\eqref{L1}$. We have that
$$
e^{\lambda_2(t,x,\xi)} \left\{
\begin{array}{ll}
\in S^{\infty}_{\mu; \frac{1}{2(1-\sigma_2)}}(\R^2), & \sigma_2 \in (\frac{1}{2}, 1),\\
\text{is of finite order}, & \sigma_2 = 1,\\
 \in S^{0}_{\mu}(\R^2), & \sigma_2 >1,
\end{array}\right.$$
whereas 
$$
e^{\lambda_1(t,x,\xi)} \left\{
\begin{array}{ll}
\in S^{\infty}_{\mu; \frac{1}{1-2\sigma_1}}(\R^2), & \sigma_2 \in (0, \frac{1}{2}),
\\
\in S^{0}_{\mu}(\R^2), & \sigma_1 \geq \frac{1}{2}.
\end{array}
\right.$$
Then, we have well-posedness in 
$$
\begin{cases}
	H^{\infty}_{\theta}, \quad \sigma_2 \in (\frac{1}{2}, 1), \sigma_1 \in (0,\frac{1}{2}), 1 < \theta < \min \{ \frac{1}{2(1-\sigma_2)}, \frac{1}{1-2\sigma_1}\},  \\
	H^{\infty}_{\theta}, \quad \sigma_2 \in (\frac{1}{2}, 1), \sigma_1 \geq \frac{1}{2},  1 < \theta < \frac{1}{2(1-\sigma_2)},  \\
	H^{\infty}_{\theta}, \quad \sigma_2 \geq 1, \sigma_1 \in (0,\frac{1}{2}), 1 < \theta < \frac{1}{1-2\sigma_1}, \\
	H^{\infty}, \quad \sigma_2 = 1, \sigma_1 \geq \frac{1}{2}, \\
	L^2, \quad \;\;\sigma_2 > 1, \sigma_1 \geq \frac{1}{2}.
\end{cases}
$$

\end{document}